%% file: waveguide.tex
\documentclass[journal]{IEEEtran}
\usepackage{cite}
\usepackage{amssymb}
\usepackage[pdftex]{graphicx}
\usepackage{amsfonts}
\usepackage{lineno}

\input latex_defs.tex

\usepackage[cmex10]{amsmath}
\usepackage{array}
\usepackage{stfloats}
\usepackage{url}
\usepackage{epstopdf}
\usepackage[square,numbers,sort&compress]{natbib}
\usepackage{subfigure}

\input latex_defs.tex
\begin{document}
%
\title{A Sufficient Condition of Having Independent TE and TM Modes in a Waveguide Filled with Homogenous Anisotropic Lossless Medium}
\author{Wei~Jiang,~\IEEEmembership{}
        Jie~Liu,~\IEEEmembership{}
        and~Qing~Huo Liu,~\IEEEmembership{Fellow,~IEEE}
\thanks{W. Jiang and J. Liu are with Institute of Electromagnetics and Acoustics, and Department of Electronic Science, Xiamen University, Xiamen 361005, China}
\thanks{Q. H. Liu is Department of Electrical and Computer Engineering,
Duke University, Durham, NC 27708 USA}}
\maketitle

\begin{abstract}
Based on the idea of the Abelian group theory in mathematics,
this paper finds a sufficient condition of having independent TE and TM modes
in a waveguide filled with homogenous anisotropic lossless medium.
For independent TE modes, we prove the nonzero cut-off wavenumbers obtained from longitudinal scalar magnetic field stimulation and transverse vector electric field stimulation are same in theory.
For independent TM modes, we also prove the nonzero cut-off wavenumbers obtained from longitudinal scalar electric field stimulation and transverse vector magnetic field stimulation are same in theory.
Finally we carry out several numerical experiments to verify the correctness of the condition given by us.
We hope that this condition is useful for the designs of waveguide with homogenous anisotropic lossless medium in microwave engineering community.
\end{abstract}

\begin{IEEEkeywords}
Waveguide Problems, Abelian Group, Independent TE Modes, Independent TM Modes. Finite Element Method.
\end{IEEEkeywords}

\IEEEpeerreviewmaketitle

\section{Introduction}
In microwave engineering, waveguide is one of very important source-free equipments,
which guides wave propagation direction. It has many applications, for example,  Hirokawa J. and Ando M. \cite{hirokawa1998} have proposed a novel feed structure to excite a plane TEM wave in a parallel-plate waveguide. The electromagnetic wave in the waveguide must be a solution of Maxwell's equations with source-free. \\
\indent According to the electromagnetic theory \cite{Balanis}, it is well-known that the waveguide filled with homogenous isotropic lossless medium has independent TE modes and TM modes. Suppose that the wave propagation direction in waveguide is $+\^z$. We also have known that for independent TE modes, we can employ the longitudinal component $h_z$ of the magnetic field or the transverse component $\ee_{t}$ of the electric field to stimulate them, and for independent TM modes, we can employ the longitudinal component $e_{z}$ of the electric field or the transverse component $\h_{t}$ of the magnetic field to stimulate them. For using longitudinal component $h_z$ of the magnetic field to stimulate independent TE modes, we need to solve an eigenvalue problem about Laplace operator with a Neumann boundary condition. This boundary condition is a natural boundary condition in finite element method (FEM), which is not enforced constraint in numerical computation. While for using longitudinal component $e_{z}$ of the electric field to stimulate independent TM modes, we need to solve an eigenvalue problem about Laplace operator with a Dirichlet boundary condition. This boundary condition is an essential boundary condition in FEM, which must be enforced constraint in numerical computation. For employing transverse component $\ee_{t}$ of the electric field to stimulate independent TE modes, we need to solve an eigenvalue problem about curl-curl operator with an essential boundary condition, which must be enforced constraint in numerical computation, while for employing transverse component $\h_{t}$ of the magnetic field to stimulate independent TM modes, we need to solve an eigenvalue problem about curl-curl operator with two natural boundary conditions,
which are all not enforced constraint in numerical computation. \\
\indent For the waveguide filled with homogenous anisotropic lossless medium, are there independent TE modes and TM modes in this waveguide? The answer is that when the medium parameters in the waveguide satisfy some condition, then there are independent TE and TM modes. This paper finds this condition (i.e., the following condition \uppercase\expandafter{\romannumeral1} and \uppercase\expandafter{\romannumeral2}) based on the idea of the Abelian group theory \cite{rotman2010advanced} (see Appendix B) in mathematics. Moreover, when the waveguide filled with homogenous anisotropic lossless medium has independent TE and TM modes, we prove that the nonzero cut-off wavenumbers obtained from using $h_{z}$ and $\ee_{t}$ to stimulate independent TE modes are same in theory, and the nonzero cut-off wavenumbers obtained from using $e_{z}$ and $\h_{t}$ to stimulate independent TM modes are also same in theory.\\
\indent When the cross section of the waveguide is not regular,
if we would like to know the first several propagable physical modes,
we usually turn to the help of numerical methods,
for example, FEM, finite difference method, and etc.
This paper will use standard linear FEM to solve the eigenvalue problems about elliptic differential operator of second-order in two-dimension,
and mixed FEM \cite{Brezzi1991} based on CT/LN edge element space and standard linear element space to solve the eigenvalue problems of curl-curl operator in two-dimension.
The basis functions in the CT/LN edge element space are constant tangential and liner normal (CT/LN).
This space is the lowest order edge element space,
which is the simplest edge element space.
This space is one of the N\'ed\'elec spaces \cite{nede}.
Note that it is careful to the selection of these two spaces in mixed FEM,
and the selection of these two spaces should be matched \cite{Brezzi1991}.
For the eigenvalue problems about curl-curl operator
(waveguide problems and resonant cavity problems just belong to this type of problem), the earlier papers on edge elements claim edge elements can entirely eliminate spurious modes. This is not correct \cite{davidson2010}. Because edge elements can only eliminate nonphysical spurious nonzero modes, however, edge elements can not eliminate nonphysical spurious zero modes because of ignoring the divergence-free condition. We have computed the resonant cavity problem with anisotropic lossless media and PEC successfully \cite{Jiang2015}. In the reference \cite{Jiang2015} we turn to the help of mixed FEM based on CT/LN edge element space and linear element space to remove all the nonphysical spurious modes. Numerical methods provided by this paper can remove all the nonphysical spurious modes (including nonphysical spurious zero modes) for the eigenvalue problems about curl-curl operator, because we have already enforced divergence-free condition in numerical computation.\\
\indent The outline of this paper is as follows.
In Section 2, we review the governing equations for waveguide problems according to the waveguide theory in electromagnetism, and find this sufficient condition of having independent TE and TM modes in a waveguide filled with homogenous anisotropic lossless medium. The proof of the conclusion that the nonzero cut-off wavenumbers obtained from using $h_{z}$ and $\ee_{t}$ to stimulate independent TE modes are same is given in Section 3. The proof of the conclusion that the nonzero cut-off wavenumbers obtained from using $e_{z}$ and $\h_{t}$ to stimulate independent TM modes are same is given in Section 4.  In Section 5, we support several FEM algorithms of PDEs from waveguide problems.  Finally we carry out some numerical experiments to verify that the condition given by us is correct, and the theory in Section 3 and 4 is also correct.
\section{Governing Equations for Waveguide Problems}
Suppose that the wave propagation direction in waveguide is $\^z$, where $\^z$ is unit vector along positive direction of $z$ axis. Let $\Gamma$ be cross section of the waveguide, $\^n$ be be the outward normal unit vector on the boundary $\partial\Gamma$ of the cross section $\Gamma$. As usual, the cross section $\Gamma$ is bounded. Because of the complexity of the cross section $\Gamma$, the boundary $\partial\Gamma$ may be not connected. In this case, there exist TEM modes. In practical applications, waveguide walls are usually made up of perfect electric conductor.\\
\indent Waveguide problem is a so called 2.5 dimensional problem, because the medium in the waveguide is two dimensional, while the electromagnetic field in the waveguide is three dimensional.  This paper only considers the waveguide problem filled with filled with homogenous anisotropic lossless medium. Suppose that the permittivity and permeability matrices are of the following form:
\begin{eqnarray}
  \d{\ep}&=&\left[
  \begin{array}{cc}
    \d{\ep}_{t} & 0 \\
    0 & \ep_{zz} \\
  \end{array}
\right],\label{ep}\\
\d{\mu}&=&\left[
  \begin{array}{cc}
    \d{\mu}_{t} & 0 \\
    0 & \mu_{zz} \\
  \end{array}
\right], \label{mu}
\end{eqnarray}
where $\d{\ep}_{t}$, $\d{\mu}_{t}$, $\ep_{zz}$ and $\mu_{zz}$ are independent of $z$. Because the medium is lossless, we have that $\d{\ep}$ and $\d{\mu}$ are two positive definite Hermitian matricies \cite{Chew1990}, i.e.,
$$\d{\ep}_{t}^{\dagger}=\d{\ep}_{t}, \quad\d{\mu}_{t}^{\dagger}=\d{\mu}_{t}, ~~\ep_{zz}>0, ~~\mu_{zz}>0,$$
where $\d{\ep}_{t}$ and $\d{\mu}_{t}$ have positive real eigenvalues, and $K^{\dagger}$ is conjugate transpose of matrix $K$.

In a usual waveguide problem, the frequency $f>0$ is given, then $\omega=2\pi f>0$. We need to solve propagation constant $jk_{z}$ and several propagable physical modes:
\begin{equation}\label{guide1}
    \E(x,y,z)=\ee(x,y)e^{-jk_{z}z},~\H(x,y,z)=\h(x,y)e^{-jk_{z}z},
\end{equation}
where $\ee(x,y)$ and $\h(x,y)$ are two three dimensional vectors only dependent with transversal coordinates $(x,y)$, while independent with longitudinal coordinate $z$. Assume that
\begin{eqnarray}
    &~&\ee(x,y)=\ee_{t}+\^ze_{z},\quad\h(x,y)=\h_{t}+\^zh_{z},\label{guide2} \\ &~&\nabla=\nabla_{t}+\^z\pd{}{z},\quad\nabla_{t}=\^x\pd{}{y}+\^x\pd{}{y}.\label{guidddde2}
\end{eqnarray}
It is clear that
$$\pd{\E}{z}=-jk_z\E,\quad\pd{\H}{z}=-jk_z\H.$$
Writing Maxwell's equations with source-free, we have:
\begin{eqnarray*}
  \curl\E&=&-j\omega\d{\mu}\H, \\
  \curl\H&=& j\omega\d{\ep}\E,\\
  \div(\d{\ep}\E)&=&0,\\
  \div(\d{\mu}\H)&=&0,
\end{eqnarray*}
where we have assumed that the time-harmonic factor is $e^{j\omega t}$. Substituting (\ref{ep}-\ref{guidddde2}) into the above Maxwell's equations, then we obtain
the following equations:
\begin{eqnarray}
  \curlt\ee_{t} &=& -j\omega\mu_{zz}\^zh_z\label{eqn31} \\
  -\^z\times\nabla_{t}e_z-jk_z\^z\times\ee_{t}&=& -j\omega\d{\mu}_{t}\h_{t} \label{eqn32}\\
  \curlt\h_{t} &=& j\omega\ep_{zz}\^ze_z\label{eqn33} \\
  -\^z\times\nabla_{t}h_z-jk_z\^z\times\h_{t}&=& j\omega\d{\ep}_{t}\ee_{t} \label{eqn34}\\
  \divt(\d{\ep}_{t}\ee_{t})&=&jk_{z}\ep_{zz}e_z\label{addeqn35}\\
  \divt(\d{\mu}_{t}\h_{t})&=&jk_{z}\mu_{zz}h_z\label{addeqn36}
\end{eqnarray}
Set $A=\left[
       \begin{array}{cc}
         0 & -1 \\
         1 & 0 \\
       \end{array}
     \right]$, then $A^2=-I_{2\times2}$,
where $I_{2\times2}$ is identity matrix of order two. We change $\^z\times$ in equations (\ref{eqn32}) and (\ref{eqn34}) into $A$, then we get
\begin{eqnarray}
\nabla_{t}e_{z}+jk_{z}\ee_{t}&=&-j\omega A\d{\mu}_{t}\h_{t},\label{eqn35}\\
\nabla_{t}h_{z}+jk_{z}\h_{t}&=&j\omega A\d{\ep}_{t}\ee_{t}.\label{eqn36}
\end{eqnarray}
From equations (\ref{eqn35}) and (\ref{eqn36}), we can obtain
\begin{eqnarray}
(-\omega^2A\d{\mu}_{t}A\d{\ep}_{t}-k_{z}^2I_{2\times2})\ee_{t}
=-jk_{z}\nabla_{t}e_{z}+j\omega A\d{\mu}_{t}\nabla_{t}h_{z}, \label{eqn37}\\
(-\omega^2A\d{\ep}_{t}A\d{\mu}_{t}-k_{z}^2I_{2\times2})\h_{t}
=-jk_{z}\nabla_{t}h_{z}-j\omega A\d{\ep}_{t}\nabla_{t}e_{z}. \label{eqn38}
\end{eqnarray}

Next we try to explore the condition of having independent TE and TM modes in a waveguide filled with homogenous anisotropic lossless medium, and the case in filled with homogenous isotropic lossless medium is a particular case of this condition.

Define set $G=\{X\in{\mathbb{C}^{2\times2}}: AX=XA\}$, then the set $G$ has very nice algebraic properties. It is very easy to prove that $(G,+)$ is an Abelian group, where $+$ is usual matrix addition in $\mathbb{C}^{2\times2}$ and $(G,\cdot)$ is also an group, where $\cdot$ is usual matrix  multiplication in $\mathbb{C}^{2\times2}$. Next we intensively consider the concrete form of the element in the group $(G,\cdot)$, let $X=\left[
       \begin{array}{cc}
         a & b \\
         c & d \\
       \end{array}
     \right]$ be in the group $(G,\cdot)$, then we have $\left[
       \begin{array}{cc}
         a & b \\
         c & d \\
       \end{array}
     \right]\left[
       \begin{array}{cc}
         0 & -1 \\
         1 & 0 \\
       \end{array}
     \right]=\left[
       \begin{array}{cc}
         0 & -1 \\
         1 & 0 \\
       \end{array}
     \right]\left[
       \begin{array}{cc}
         a & b \\
         c & d \\
       \end{array}
     \right]$, therefore we can get $a=d, b=-c$. When $\d{\ep}_{t}\in{G}$ and $\d{\mu}_{t}\in{G}$, then we can assume that $\d{\ep}_{t}=\left[
       \begin{array}{cc}
         \ep & aj \\
         -aj & \ep \\
       \end{array}
     \right]$, $\d{\mu}_{t}=\left[
       \begin{array}{cc}
         \mu & bj \\
         -bj & \mu \\
       \end{array}
     \right]$,
where $\ep$, $\mu$, $a$, $b$ are all real numbers, and $\d{\ep}_{t}$ and $\d{\mu}_{t}$ are two positive definite Hermitian matrices.
$$\left[
       \begin{array}{cc}
         \ep & aj \\
         -aj & \ep \\
       \end{array}
     \right]\left[
       \begin{array}{cc}
         \mu & bj \\
         -bj & \mu \\
       \end{array}
     \right]=\left[
       \begin{array}{cc}
         \ep\mu+ab & (b\ep+a\mu)j \\
         -(b\ep+a\mu)j & \ep\mu+ab \\
       \end{array}
     \right].$$
$$\left[
       \begin{array}{cc}
         \mu & bj \\
         -bj & \mu \\
       \end{array}
     \right]\left[
       \begin{array}{cc}
         \ep & aj \\
         -aj & \ep \\
       \end{array}
     \right]=\left[
       \begin{array}{cc}
         \ep\mu+ab & (b\ep+a\mu)j \\
         -(b\ep+a\mu)j & \ep\mu+ab \\
       \end{array}
     \right].$$
Then $\d{\mu}_{t}\d{\ep}_{t}=\d{\ep}_{t}\d{\mu}_{t}$ is valid, therefore $(G,\cdot)$ is also an Abelian group. Particularly, when $b\ep+a\mu=0$, then we have
\begin{gather}
    \d{\ep}_{t}\d{\mu}_{t}=\d{\mu}_{t}\d{\ep}_{t}= (\ep\mu+ab)I_{2\times2}\neq 0. \label{matrix3w1}\\
 \d{\ep}_{t}^{-1}=\frac{1}{\ep\mu+ab}\d{\mu}_{t}\label{matrix3w2}.\\
 \d{\mu}_{t}^{-1}=\frac{1}{\ep\mu+ab}\d{\ep}_{t}\label{matrix3w2}.
\end{gather}
Suppose that $\d{\ep}_{t}\in{G}$, $\d{\mu}_{t}\in{G}$, and they satisfies the condition $b\ep+a\mu=0$, then
$$-\omega^2A\d{\mu}_{t}A\d{\ep}_{t}
=\omega^2\d{\mu}_{t}\d{\ep}_{t}=\omega^2\d{\ep}_{t}\d{\mu}_{t}=\omega^2(\ep\mu+ab)I_{2\times2};$$
$$-\omega^2A\d{\ep}_{t}A\d{\mu}_{t}
=\omega^2\d{\ep}_{t}\d{\mu}_{t}=\omega^2\d{\mu}_{t}\d{\ep}_{t}=\omega^2(\ep\mu+ab)I_{2\times2}.$$
Thus coefficient matrices in equations (\ref{eqn37}) and (\ref{eqn38}) are same:
$$-\omega^2A\d{\mu}_{t}A\d{\ep}_{t}-k_{z}^2I_{2\times2}=(\omega^2\ep\mu+\omega^2ab-k_{z}^2)I_{2\times2};$$
$$-\omega^2A\d{\ep}_{t}A\d{\mu}_{t}-k_{z}^2I_{2\times2}=(\omega^2\ep\mu+\omega^2ab-k_{z}^2)I_{2\times2}.$$
Define
\begin{eqnarray}
k_{t}^2=k^2-k_{z}^2, \quad k^2=\omega^2(\ep\mu+ab),\label{3www1}
\end{eqnarray}
where $k_t$ is called cut-off wavenumber, and $k_z$ is called phase constant.
From equations (\ref{eqn37})-(\ref{eqn38}), when $k_{t}\neq0$, we can infer that
\begin{eqnarray}
\ee_{t}
=-\frac{jk_{z}}{k_{t}^2}\nabla_{t}e_{z}+\frac{j\omega}{k_{t}^2}\^z\times(\d{\mu}_{t}\nabla_{t}h_{z}), \label{eqn39}\\
\h_{t}
=-\frac{jk_{z}}{k_{t}^2}\nabla_{t}h_{z}-\frac{j\omega}{k_{t}^2}\^z\times(\d{\ep}_{t}\nabla_{t}e_{z}). \label{eqn310}
\end{eqnarray}
From (\ref{eqn39}) and (\ref{eqn31}), we can get:
\begin{eqnarray}
-\divt(\d{\mu}_{t}\nabla_{t}h_{z})=k_{t}^2\mu_{zz}h_z.\label{eqn41}
\end{eqnarray}
From (\ref{eqn310}) and (\ref{eqn33}), we can get:
\begin{eqnarray}
-\divt(\d{\ep}_{t}\nabla_{t}e_{z})=k_{t}^2\ep_{zz}e_z.\label{eqn42}
\end{eqnarray}
if $e_{z}=0$, $h_{z}=0$ and equations (\ref{eqn31})-(\ref{addeqn36}) have nontrivial solution, then we must have $k_{t}^2=k^2-k_{z}^2=0$, thus, $k_{z}=\omega\sqrt{\ep\mu+ab}~(\ep\mu+ab>0)$. The transverse components of electromagnetic field are directly obtained from the equations (\ref{addeqn35})-(\ref{addeqn36}). This is so called TEM mode.

In conclusion, when the medium parameters in a waveguide satisfies the following condition:
\begin{enumerate}
  \item [\uppercase\expandafter{\romannumeral1}]\begin{eqnarray*}
            \d{\ep}=\left[
  \begin{array}{ccc}
    {\ep} &aj &0 \\
    -aj & \ep&0 \\
    0 & 0&\ep_{zz} \\
  \end{array}
\right],
\quad\d{\mu}=\left[
  \begin{array}{ccc}
     \mu &bj& 0 \\
     -bj&\mu& 0 \\
    0 & 0&\mu_{zz} \\
  \end{array}
\right]
\end{eqnarray*}
are two positive definite Hermitian matrices, and they are constant matrices.
  \item [\uppercase\expandafter{\romannumeral2}]$$b\ep+a\mu=0,$$
\end{enumerate}
 then the independent TE and TM modes must exist in the waveguide filled with homogenous anisotropic lossless medium.\\
\noindent{\bf{Remark:}} It is well-known that when the waveguide is filled with homogenous isotropic lossless medium,
then there are independent TE and TM modes in this waveguide.
In this case, $a=0$, $b=0$, $\ep=\ep_{zz}>0$ and $\mu=\mu_{zz}>0$.
Obviously, the medium parameters have already satisfied the above condition \uppercase\expandafter{\romannumeral1} and \uppercase\expandafter{\romannumeral2}.

Next we consider the boundary condition about waveguide problem.
Because the waveguide walls are usually made up of metal,
the boundary condition is usually perfect electric conductor (PEC) boundary condition. From electromagnetic theory \cite{Balanis}, we have $$\^n\times\E=0,\quad\^n\cdot\B=0.$$

For the TE modes, if we use longitudinal scalar field $h_{z}$ to stimulate them, then the boundary condition will be $\^n\cdot(\d{\mu}_{t}\gradt h_z)=0\mbox{~on~}\partial\Gamma$;
if we use transverse vector field $\ee_{t}$ to stimulate them,
then the boundary condition will be $\^n\times\ee_{t}=0\mbox{~on~}\partial\Gamma.$

For the TM modes, if we use longitudinal scalar field $e_{z}$ to stimulate them, then the boundary condition will be $e_{z}=0\mbox{~on~}\partial\Gamma$;
if we use transverse vector field $\h_{t}$ to stimulate them, then the boundary condition will be $\^n\cdot(\d{\mu}_{t}\h_{t})=0\mbox{~on~}\partial\Gamma$ and $\^n\times(\ep_{zz}^{-1}\curlt\h_{t})=0~~\mbox{on}~\partial\Gamma.$
\section{Independent TE Modes}
In this section, the medium parameters in the waveguide we consider satisfy the condition \uppercase\expandafter{\romannumeral1} and \uppercase\expandafter{\romannumeral2}.
\subsection{Simulate TE Modes using Longitudinal Component $h_z$}
\indent Seek $k_{t}\in{\mathbb{R}}$, $h_z\neq 0$, such that
\begin{equation} \label{3eqs1}
\left\{ \begin{aligned}
         -\divt(\d{\mu}_{t}\nabla_{t}h_{z})&=k_{t}^2\mu_{zz}h_z~~\mbox{in}~\Gamma, \\
         \^n\cdot(\d{\mu}_{t}\gradt h_z)&=0~~\mbox{on}~\partial\Gamma.
        \end{aligned} \right.
\end{equation}
The PDE (\ref{3eqs1}) is an eigenvalue problem about elliptic partial differential operator of second order in two-dimension,
and this problem is with Neumann boundary condition. The researches about this problem are very thorough in computational mathematics. For details, please see \cite{Boffi1}.\\
\indent It is clear that the equation (\ref{3eqs1}) has a particular solution $k_t=0$, $h_z=1$, but this solution is not a propagable physical mode, and this is a spurious mode. We shall not consider this solution. Once $(k_{t},h_z)$ is solved, then
\begin{eqnarray*}
\ee_{t}&=&\frac{j\omega}{k_{t}^2}\^z\times(\d{\mu}_{t}\nabla_{t}h_{z});\\
\h_{t}&=&-\frac{jk_{z}}{k_{t}^2}\nabla_{t}h_{z}.
\end{eqnarray*}
From the spectral theory of self-adjoint compact operator \cite{chatelin1983}, we have already known that PDE (\ref{3eqs1}) has countable real eigenvalues. i.e.,
$$0=k_{t}^{(0)}<k_{t}^{(1)}\leq k_{t}^{(2)}\leq\cdots\leq\cdots\rightarrow+\infty.$$
Denote $k_{t}^{(small)}=k_{t}^{(1)}$.
Because $k_{z}^2=k^2-k_{t}^2$, when $k>k_{t}^{(small)}$, i.e., $f>f_{t}^{(small)}$,
where $f_{t}^{(small)}$ is the smallest cut-off frequency in the waveguide,
then the waveguide has propagable physical modes. Therefore the waveguide has the function of high-pass filter. That is the electromagnetic wave with low-frequency can not propagate in the waveguide, while the electromagnetic wave with high-frequency can propagate in the waveguide.\\
\noindent{\bf{Lemma 1.}} If $k_{t}\neq0$ and $(k_{t},h_{z})$ is an eigen-pair of PDE (\ref{3eqs1}),  then $(k_{t}, \frac{j\omega}{k_{t}^2}\^z\times(\d{\mu}_{t}\nabla_{t}h_{z}))$ is also an eigen-pair of PDEs (\ref{3eqs2}). \\
\noindent Proof. Please see Appendix A.
\subsection{Simulate TE Modes using Transverse Component $\ee_t$}
\indent Find $k_{t}\in{\mathbb{R}}$, $0\neq\ee_t$, such that
\begin{equation} \label{3eqs2}
\left\{ \begin{aligned}
         \curlt\bigg(\mu_{zz}^{-1}\curlt\ee_t\bigg)&=k_{t}^2\d{\mu}_{t}^{-1}\ee_t~~\mbox{in}~\Gamma, \\
         \divt(\d{\mu}_{t}^{-1}\ee_t)&=0~~\mbox{in}~\Gamma,\\
         \^n\times\ee_t&=0~~\mbox{on}~\partial\Gamma.
        \end{aligned} \right.
\end{equation}
The PDEs (\ref{3eqs2}) is an eigenvalue problem about curl-curl operator in two-dimension.
Obviously, when $k_t\neq0$, the second equation of PDEs (\ref{3eqs2}) can be derived from the first equation of PDEs (\ref{3eqs2}),
this is because $0=\divt\bigg(\curlt\mu_{zz}^{-1}\curlt\ee_t\bigg)=k_{t}^2\divt(\d{\mu}_{t}^{-1}\ee_t)$,
then we have $\divt(\d{\mu}_{t}^{-1}\ee_t)=0$.
But we can not omit the divergence-free condition $\divt(\d{\mu}_{t}^{-1}\ee_t)=0$ in numerical computation,
otherwise PDEs (\ref{3eqs2}) will introduce spurious zero modes \cite{jiang2014}. Because $\divt(\d{\mu}_{t}^{-1}\ee_t)=0$ and (\ref{matrix3w2}), we can get $\frac{1}{\ep\mu+ab}\divt(\d{\ep}_{t}\ee_t)=0$, then $\divt(\d{\ep}_{t}\ee_t)=0$, which is just Gauss' law for electric field. In fact, $\divt(\d{\mu}_{t}^{-1}\ee_t)=0$ and $\divt(\d{\ep}_{t}\ee_t)=0$ are equivalent.

When the boundary $\partial\Gamma$ is not connected, then there are TEM modes in this waveguide. Suppose that $\partial\Gamma=\bigcup_{k=1}^{N}\partial\Gamma_{k}$, $\partial\Gamma_{i}\bigcap\partial\Gamma_{j}=\varnothing(i\neq j)$, where $\partial\Gamma_{k} (k=1,2,\cdots,N)$ is the k-th connected boundary. Then number of TEM modes is $N-1$. Which can be verified by the following numerical experiments. PDEs (\ref{3eqs2}) can stimulate these TEM modes, while PDE (\ref{3eqs1}) can not stimulate these TEM modes, which is an advantage of using $\ee_{t}$ to simulate waveguide modes. Moreover, PDE (\ref{3eqs1}) will also introduce a nonphysical spurious zero mode.

Once $(k_{t},\ee_{t})$ is solved, then
$$\h_{t}=\frac{k_{z}\d{\mu}_{t}^{-1}(\^z\times\ee_t)}{\omega}=\frac{k_{z}\^z\times(\d{\mu}_{t}^{-1}\ee_t)}{\omega};$$
$$h_z=\^z\cdot\frac{j\curlt\ee_t}{\omega\mu_{zz}}.$$

\noindent{\bf{Lemma 2.}} If $k_{t}\neq0$ and $(k_{t},\ee_{t})$ is an eigen-pair of PDEs (\ref{3eqs2}),  then $(k_{t}, \^z\cdot\frac{j\curlt\ee_t}{\omega\mu_{zz}})$ is also an eigen-pair of PDE (\ref{3eqs1}). \\
\noindent Proof. Please see Appendix A.\\
\indent In a word, we have already obtained the following theorem. \\
\noindent{\bf{Theorem 1.}} If the medium parameters in a waveguide satisfies the condition \uppercase\expandafter{\romannumeral1} and \uppercase\expandafter{\romannumeral2}, then there exist independent TE modes in this waveguide, and nonzero cut-off wavenumber $k_{t}$ between PDE (\ref{3eqs1}) and PDEs (\ref{3eqs2}) is same in theory.

\section{Independent TM Modes}
In this section, the medium parameters in the waveguide we consider also have the condition \uppercase\expandafter{\romannumeral1} and \uppercase\expandafter{\romannumeral2}.
\subsection{Simulate TM Modes using Longitudinal Component $e_z$}
\indent Seek $k_{t}\in{\mathbb{R}}$, $e_z\neq 0$, such that
\begin{equation} \label{3eqsm1}
\left\{ \begin{aligned}
         -\divt(\d{\ep}_{t}\nabla_{t}e_{z})&=k_{t}^2\ep_{zz}e_z~~\mbox{in}~\Gamma, \\
         e_z&=0~~\mbox{on}~\partial\Gamma.
        \end{aligned} \right.
\end{equation}
The PDE (\ref{3eqsm1}) is an eigenvalue problem about elliptic partial differential operator of second order in two-dimension, and this problem is with Dirichlet boundary condition.

\indent It is easy to prove that $k_{t}=0$ is not an eigenvalue of PDE (\ref{3eqsm1}).
In fact, if $k_{t}=0$ is an eigenvalue of PDE (\ref{3eqsm1}),
then we can prove that $e_{z}=0\mbox{~in~}\Gamma$ by virtue of the Poincar\'{e} inequality $\|e_{z}\|_{1}\leq C\|\grad e_{z}\|_{0}$ \cite{ciar}.
According to spectral theory of self-adjoint compact operator \cite{chatelin1983}, the PDE (\ref{3eqsm1}) has countable real positive eigenvalues.
$$0<k_{t}^{(1)}\leq k_{t}^{(2)}\leq k_{t}^{(3)}\leq\cdots\leq\cdots\rightarrow+\infty.$$
 Once $(k_{t},e_{z})$ is solved, then
 \begin{eqnarray*}
\ee_{t}
&=&-\frac{jk_{z}}{k_{t}^2}\nabla_{t}e_{z};\\
\h_{t}
&=&-\frac{j\omega}{k_{t}^2}\^z\times(\d{\ep}_{t}\nabla_{t}e_{z}).
\end{eqnarray*}
\noindent{\bf{Lemma 3.}} If $k_{t}\neq0$ and $(k_{t},e_{z})$ is an eigen-pair of PDE (\ref{3eqsm1}),  then $(k_{t}, -\frac{j\omega}{k_{t}^2}\^z\times(\d{\ep}_{t}\nabla_{t}e_{z}))$ is also an eigen-pair of PDEs (\ref{3eqsm2}). \\
\indent Proof. For the verification of the fist two equations in PDEs (\ref{3eqsm2}), this is the same as the TE modes, therefore we omit this step. We mainly need to verify the correctness of boundary conditions in PDEs (\ref{3eqsm2}). According to $\^n\times\E=0$ and $e_z=0$ on $\partial\Gamma$, then we have $\^n\times(\ee_{t}+\^ze_{z})e^{-jk_{z}z}=\^n\times\ee_{t}e^{-jk_{z}z}=0$ on $\partial\Gamma$. Therefore we have $\^n\times\ee_{t}=0$ on $\partial\Gamma$. For TM modes, we have $h_{z}=0$. We can take $h_{z}=0$ in the equation (\ref{eqn34}), then we obtain
$$\ee_{t}=-\frac{k_{z}\d{\ep}_{t}^{-1}(\^z\times\h_t)}{\omega}=-\frac{k_{z}\^z\times(\d{\ep}_{t}^{-1}\h_t)}{\omega}.$$ Then we have
\begin{eqnarray*}
  0 &=& \^n\times\ee_{t}=-\frac{k_{z}}{\omega}\^n\times\bigg(\^z\times(\d{\ep}_{t}^{-1}\h_t)\bigg)\\
  &=&-\frac{k_{z}}{\omega}\bigg[\^z(\^n\cdot(\d{\ep}_{t}^{-1}\h_t))-\d{\ep}_{t}^{-1}\h_t(\^n\cdot\^z)\bigg]\\
  &=&-\frac{k_{z}}{\omega}\^z(\^n\cdot(\d{\ep}_{t}^{-1}\h_t)) \mbox{~on~}\partial\Gamma.
\end{eqnarray*}
Therefore $\^n\cdot(\d{\ep}_{t}^{-1}\h_t)=0\mbox{~on~}\partial\Gamma.$
In addition, according to (\ref{eqn33}) and $e_z=0$ on $\partial\Gamma$, we have
\begin{eqnarray*}
  \^n\times(\ep_{zz}^{-1}\curlt\h_{t})=\^n\times(j\omega\^ze_{z})=j\omega\^n\times\^ze_{z}=0 \mbox{~on~} \partial\Gamma.
\end{eqnarray*}
The proof of {\bf Lemma 3} is completed.
\subsection{Simulate TM Modes using Transverse Component $\h_t$}
\indent Find $k_{t}\in{\mathbb{R}}$, $0\neq\ee_t$, such that
\begin{equation} \label{3eqsm2}
\left\{ \begin{aligned}
         \curlt\bigg(\ep_{zz}^{-1}\curlt\h_t\bigg)&=k_{t}^2\d{\ep}_{t}^{-1}\h_t~~\mbox{in}~\Gamma, \\
         \divt(\d{\ep}_{t}^{-1}\h_t)&=0~~\mbox{in}~\Gamma, \\
         \^n\times(\ep_{zz}^{-1}\curlt\h_{t})&=0~~\mbox{on}~\partial\Gamma,\\
         \^n\cdot(\d{\ep}_{t}^{-1}\h_{t})&=0~~\mbox{on}~\partial\Gamma.
        \end{aligned} \right.
\end{equation}
The PDEs (\ref{3eqsm2}) is an eigenvalue problem about curl-curl operator in two-dimension.
Note that $\divt(\d{\ep}_{t}^{-1}\h_t)=0$ and $\divt(\d{\mu}_{t}\h_t)=0$ are equivalent because of the equation (\ref{matrix3w2}), and $\^n\cdot(\d{\ep}_{t}^{-1}\h_{t})=0$ and
$\^n\cdot(\d{\mu}_{t}\h_{t})=0$ are also equivalent because of the equation (\ref{matrix3w2}).\\
 \indent When $\partial\Gamma$ is not connected, the PDEs (\ref{3eqsm2}) can also stimulate TEM modes, while PDE (\ref{3eqsm1}) can not stimulate these TEM modes, which is an advantage of using $\h_{t}$ to simulate waveguide modes.

Once $(k_{t},\h_{t})$ is solved, then
\begin{eqnarray}
  \ee_{t}&=&-\frac{k_{z}\d{\ep}_{t}^{-1}(\^z\times\h_t)}{\omega}=-\frac{k_{z}\^z\times(\d{\ep}_{t}^{-1}\h_t)}{\omega}; \label{boundary4}\\
\^ze_z&=&\frac{\ep_{zz}^{-1}\curlt\h_t}{j\omega}\label{boundary45}.
\end{eqnarray}
\noindent{\bf{Lemma 4.}} If $k_{t}\neq0$ and $(k_{t},\h_{t})$ is an eigen-pair of PDEs (\ref{3eqsm2}),  then $(k_{t}, \^z\cdot\frac{\curlt\h_t}{j\omega\ep_{zz}})$ is also an eigen-pair of PDE (\ref{3eqsm1}). \\
\indent Proof. For the verification of the equation in PDE (\ref{3eqsm1}),
this is the same as the TE modes,
therefore we omit this step.
We mainly need to verify the correctness of boundary condition in PDE (\ref{3eqsm1}).
From (\ref{boundary4}) and the boundary conditions in PDEs (\ref{3eqsm2}), we can get
\begin{eqnarray*}
  &&\^n\times\ee_{t}=\^n\times(-\frac{k_{z}\^z\times(\d{\ep}_{t}^{-1}\h_t)}{\omega}) \\
  &=&-\frac{k_{z}}{\omega}\^n\times(\^z\times(\d{\ep}_{t}^{-1}\h_t))\\
  &=&\frac{k_{z}}{\omega}\bigg(\d{\ep}_{t}^{-1}\h_t(\^n\cdot\^z)-\^z(\^n\cdot(\d{\ep}_{t}^{-1}\h_t))\bigg)\\
  &=&0 \mbox{~on~}\partial\Gamma.
\end{eqnarray*}
According to $\^n\times\E=0$ on $\partial\Gamma$, then we have $\^n\times(\ee_{t}+\^ze_{z})e^{-jk_{z}z}=0$ on $\partial\Gamma$. From the above equation, then we get $\^n\times(\^ze_{z})=0$ on $\partial\Gamma$, thus $e_z=0$ on $\partial\Gamma$. In addition, according to (\ref{boundary45}), we have
\begin{eqnarray*}
 \^n\times\^ze_z =\frac{1}{j\omega}\^n\times(\ep_{zz}^{-1}\curlt\h_{t})=0\mbox{~on~}\partial\Gamma,
\end{eqnarray*}
therefore we prove $e_z=0\mbox{~on~}\partial\Gamma$ again. The proof of {\bf Lemma 4} is completed.

In a word, we have already obtained the following theorem. \\
\noindent{\bf{Theorem 2.}} If the medium parameters in a waveguide satisfy the condition \uppercase\expandafter{\romannumeral1} and \uppercase\expandafter{\romannumeral2}, then there exist independent TM modes in this waveguide, and nonzero cut-off wavenumber $k_{t}$ between PDE (\ref{3eqsm1}) and PDEs (\ref{3eqsm2}) is same in theory.
\section{Numerical Methods for Independent TE Modes and TM Modes}
\subsection{Numerical Method for the scalar PDE (\ref{3eqs1}) and (\ref{3eqsm1})}
Firstly, we consider the numerical treatments for the scalar PDE (\ref{3eqs1}) and (\ref{3eqsm1}). There are many finite element methods to deal with the scalar PDE (\ref{3eqs1}) and (\ref{3eqsm1}) now. For example, conforming FEM, nonconforming FEM, and mixed FEM can solve all them numerically. FEM is a numerical method based on the variational form of PDE. Thus we need to give the corresponding variational forms about the scalar PDE (\ref{3eqs1}) and (\ref{3eqsm1}) firstly.

As usual, we need to introduce some Hilbert spaces on complex field $\mathbb{C}$ associated with PDE (\ref{3eqs1}) and (\ref{3eqsm1}):
\begin{gather*}
     L^{2}(\Gamma)=\big\{f: \int_{\Gamma}|f(x,y)|^2dxdy<+\infty\big\},\\
    H^{1}(\Gamma)=\big\{v\in{L^2(\Gamma):\gradt v\in{(L^2(\Gamma))^2}}\big\},\\
    H_{0}^1(\Gamma)=\big\{v\in{H^1(\Gamma)}:v|_{\p\Gamma}=0\big\}.
\end{gather*}
The standard inner products and norms in the above Hilbert spaces are defined as following:
\begin{gather*}
(u,v)_{0}=\int_{\Gamma}{u\overline{v}dxdy},~~\forall~u,v\in{L^2(\Gamma)},\\
\|u\|_{0}=\sqrt{(u,v)_{0}}=\big(\int_{\Gamma}{|u(x,y)|^2dxdy}\big)^{\frac{1}{2}},~~\forall~u\in{L^2(\Gamma)},\\
(u,v)_{1}=\int_{\Gamma}{(u\overline{v}+\gradt u\cdot\gradt \overline{v}) dxdy},~~\forall~u,v\in{H^1(\Gamma)},\\
\|v\|_{1}=\sqrt{(u,v)_{1}}=\big(\|v\|_{0}^2+\|\gradt v\|_{0}^2\big)^{\frac{1}{2}},~~\forall~v\in{H^1(\Gamma)},
\end{gather*}
where the function $\overline{v}$ is complex conjugation of function ${v}$.
\indent Denote the continuous sesquilinear forms:
\begin{eqnarray*}
a_{1}:&&H^{1}(\Gamma)\times H^{1}(\Gamma)\rightarrow{\mathbb{C}}\\
&&(h_z,v)\rightarrow\int_{\Gamma}{\d{\mu}_{t}\gradt h_z\cdot\gradt \overline{v}}dxdy,\\
a_{2}:&&H_{0}^{1}(\Gamma)\times H_{0}^{1}(\Gamma)\rightarrow{\mathbb{C}}\\
&&(e_z,v)\rightarrow\int_{\Gamma}{\d{\ep}_{t}\gradt e_z\cdot\gradt \overline{v}}dxdy,\\
b_{1}:&& L^2(\Gamma)\times L^2(\Gamma)\rightarrow{\mathbb{C}}\\
&& (h_z,v)\rightarrow\int_{\Gamma}{\mu_{zz} h_z \overline{v}}dxdy,\\
b_{2}:&& L^2(\Gamma)\times L^2(\Gamma)\rightarrow{\mathbb{C}}\\
&& (e_z,v)\rightarrow\int_{\Gamma}{\ep_{zz} e_z \overline{v}}dxdy.\\
\end{eqnarray*}

Using scalar Green's formula (\ref{scalargreen1}), we get the variational forms of PDE (\ref{3eqs1}) and (\ref{3eqsm1}), respectively.\\
\indent For the TE modes, find $k_{t}\in{\mathbb{R}}$, $0\neq h_z\in{H^{1}(\Gamma)}$, such that
\begin{gather}\label{var1}
  a_{1}(h_z,v)=k_{t}^2 b_{1}(h_z,v),~~\forall~v\in{H^{1}(\Gamma)}.
\end{gather}
\indent For the TM modes, find $k_{t}\in{\mathbb{R}}$, $0\neq e_z\in{H_{0}^{1}(\Gamma)}$, such that
\begin{gather}\label{var2}
  a_{2}(e_z,v)=k_{t}^2 b_{2}(e_z,v),~~\forall~v\in{H_{0}^{1}(\Gamma)}.
\end{gather}
\indent Let $\mathcal{T}_{h}$ be a regular triangular partition \cite{ciar} of $\Gamma$ with mesh parameter $h$, where $h$ is the measure of the mesh intensive. The advantage of triangular mesh is that it can approximate arbitrary bounded domain $\Gamma$ in two-dimension well. We are only plan to use the standard linear element space $V_{1}=S^{h}$ to approximate the Hilbert space $H^{1}(\Gamma)$, and $V_{2}=S^{h}\bigcap H_{0}^{1}(\Gamma)$ to approximate the Hilbert space $H_{0}^{1}(\Gamma)$. It is clear that $V_1$ is a FEM subspace of $H^{1}(\Gamma)$, and $V_2$ is a FEM subspace of $H_{0}^{1}(\Gamma)$. Restricting (\ref{var1}) and (\ref{var2}) on $V_{1}$ and $V_{2}$  respectively, then we get discrete variational forms about (\ref{var1}) and (\ref{var2}).\\
\indent For the TE modes, seek $k_{t,h}\in{\mathbb{R}}$, $0\neq h_{z}^{h}\in{V_{1}}$, such that
\begin{gather}\label{disvar1}
  a_{1}(h_{z}^{h},v)=k_{t,h}^{2} b_{1}(h_{z}^{h},v),~~\forall~v\in{V_{1}}.
\end{gather}
\indent For the TM modes, seek $k_{t,h}\in{\mathbb{R}}$, $0\neq e_{z}^{h}\in{V_{2}}$, such that
\begin{gather}\label{disvar2}
  a_{2}(e_{z}^{h},v)=k_{t,h}^2 b_{2}(e_{z}^{h},v),~~\forall~v\in{V_{2}}.
\end{gather}
where superscript or subscript $h$ stands for approximate solution of the original problem in order to distinguish exact solution of the original problem.
The above method is called the conforming FEM.
Assume that the approximate eigenvalues of variational form (\ref{disvar1}) are
$$0\approx k_{t,h}^{(0)}<k_{t,h}^{(1)}\leq k_{t,h}^{(2)}\leq\cdots\leq k_{t,h}^{(n-1)},$$
where $n$ is all the nodal number in the mesh $\mathcal{T}_{h}$, including the number of the nodes on the boundary $\partial\Gamma$. According to min-max principle \cite{strang1973}, we have
$$k_{t}^{(l)}\leq k_{t,h}^{(l)},~l=1,2,\cdots,n-1.$$
Similarly, suppose that the approximate eigenvalues of variational form (\ref{disvar2}) are
$$0<k_{t,h}^{(1)}\leq k_{t,h}^{(2)}\leq\cdots\leq k_{t,h}^{(m)},$$
where $m$ is interior nodal number in the mesh $\mathcal{T}_{h}$, excluding the number of the nodes on the boundary $\partial\Gamma$. According to min-max principle \cite{strang1973}, we have
$$k_{t}^{(l)}\leq k_{t,h}^{(l)},~l=1,2,\cdots,m.$$
Therefore the numerical eigenvalues obtained from conforming FEM approximate exact eigenvalues from above, which is the characteristic of employing conforming FEM to solve eigenvalue problems about the elliptic operator of second order. From numerical experiments, we shall find that all the numerical eigenvalues $k_{t,h}^{(l)},~l=1,2,3,\cdots$ are monotone decreasing as the decrease of mesh parameter $h$. \\
\indent Let $\{\phi_{k}\}_{k=1}^{n}$ be basis functions on the space $V_{1}$ and $h_{z}^{h}=\sum_{k=1}^{n}\xi_{k}\phi_{k}$. As last, we need to solve the following generalized matrix eigenvalue problem:
\begin{equation}\label{TEscalar}
    A_{1}\xi=k_{t,h}^2B_{1}\xi,
\end{equation}
where
\begin{gather*}
  A_{1}=(a_{ij}^{(1)})\in\mathbb{C}^{n\times n},~~B_{1}=(b_{ij}^{(1)})\in\mathbb{C}^{n\times n};\\
  a_{ij}^{(1)}=a_{1}(\phi_{j},\phi_{i}),~~b_{ij}^{(1)}=b_{1}(\phi_{j},\phi_{i});\\
  \xi=[\xi_{1},\xi_{2},\cdots,\xi_{n}]^{T}.
\end{gather*}
\indent Let $\{\varphi_{k}\}_{k=1}^{m}$ be basis functions on the space $V_{2}$ and $e_{z}^{h}=\sum_{k=1}^{m}\xi_{k}\varphi_{k}$. As last, we need to solve the following generalized matrix eigenvalue problem:
\begin{equation}\label{TMscalar}
    A_{2}\xi=k_{t,h}^2B_{2}\xi,
\end{equation}
where
\begin{gather*}
  A_{2}=(a_{ij}^{(2)})\in\mathbb{C}^{m\times m},~~B_{2}=(b_{ij}^{(2)})\in\mathbb{C}^{m\times m};\\
  a_{ij}^{(2)}=a_{2}(\varphi_{j},\varphi_{i}),~~b_{ij}^{(2)}=b_{2}(\varphi_{j},\varphi_{i});\\
  \xi=[\xi_{1},\xi_{2},\cdots,\xi_{m}]^{T}.
\end{gather*}
\subsection{Numerical Method for the Vector PDEs (\ref{3eqs2}) and (\ref{3eqsm2})}
Secondly, we consider the numerical treatments for the vector PDEs (\ref{3eqs2}) and (\ref{3eqsm2}), which are much more difficult than for the scalar PDE (\ref{3eqs1}) and
(\ref{3eqsm1}), because of the constraint of divergence-free condition
in PDEs (\ref{3eqs2}) and (\ref{3eqsm2}).
Until now the conforming FEM discretization about
PDEs (\ref{3eqs2}) and (\ref{3eqsm2}) has not appear,
because it is very difficult to construct a conforming
FEM subspace corresponding to the vector PDEs (\ref{3eqs2}) and (\ref{3eqsm2}) \cite{boffi1999}. We shall use mixed FEM to solve the vector PDEs (\ref{3eqs2}) and (\ref{3eqsm2}). Next we give mixed variational forms associated with the vector PDEs (\ref{3eqs2}) and (\ref{3eqsm2}).\\
\indent As usual, we need to introduce some Hilbert spaces on complex field $\mathbb{C}$ associated with PDEs (\ref{3eqs2}) and (\ref{3eqsm2}):
\begin{gather*}
    \H(\mbox{curlt},\Gamma)=\big\{\F\in{(L^2(\Gamma))^2}: \curlt{\F}\in{L^2(\Gamma)}\big\},\\
    \H_{0}(\mbox{curlt},\Gamma)=\big\{\F\in{\H(\mbox{curlt},\Gamma)}: \^n\times{\F}|_{\p\Gamma}=0\big\}.
\end{gather*}
The standard inner product and norm in the above Hilbert space is defined as following: for $\forall~\F_{1},\F_{2},\F\in{\H(\mbox{curlt},\Gamma)}$,
\begin{gather*}
(\F_{1},\F_{2})_{\mbox{curlt}}=\int_{\Gamma}{\F_{1}\cdot\overline{\F_{2}}+\curlt\F_{1}\cdot\curlt\overline{\F_{2}}dxdy},\\
\|\F\|_{\mbox{curlt}}=\sqrt{(\F,\F)_{\mbox{curlt}}}=\big(\|\F\|_{0}^2+\|\curlt\F\|_{0}^2\big)^{\frac{1}{2}}.
\end{gather*}
Denote the continuous sesquilinear forms:
\begin{eqnarray*}
\mathcal{A}_{1}:&&\H_{0}(\mbox{curlt},\Gamma)\times\H_{0}(\mbox{curlt},\Gamma)\rightarrow{\mathbb{C}}\\
&&(\ee_{t},\F)\rightarrow\int_{\Gamma}{{\mu}_{zz}^{-1}\curlt\ee_{t}\cdot{\curlt\overline{\F}}}dxdy,\\
\mathcal{A}_{2}:&&\H(\mbox{curlt},\Gamma)\times\H(\mbox{curlt},\Gamma)\rightarrow{\mathbb{C}}\\
&&(\h_{t},\F)\rightarrow\int_{\Gamma}{{\ep}_{zz}^{-1}\curlt\h_{t}\cdot{\curlt\overline{\F}}}dxdy,\\
\mathcal{B}_{1}:&&(L^{2}(\Gamma))^2\times(L^{2}(\Gamma))^2\rightarrow{\mathbb{C}}\\
&&(\ee_{t},\F)\rightarrow\int_{\Gamma}{\d{\mu}_{t}^{-1}\ee_{t}\cdot{\overline{\F}}}dxdy,\\
\mathcal{B}_{2}:&&(L^{2}(\Gamma))^2\times(L^{2}(\Gamma))^2\rightarrow{\mathbb{C}}\\
&&(\h_{t},\F)\rightarrow\int_{\Gamma}{\d{\ep}_{t}^{-1}\h_{t}\cdot{\overline{\F}}}dxdy,\\
\mathcal{C}_{1}:&&\H_{0}(\mbox{curlt},\Gamma)\times{H_{0}^1(\Gamma)}\rightarrow{\mathbb{C}}\\
&&(\ee_{t},q)\rightarrow\int_{\Gamma}{\d{\mu}_{t}^{-1}\ee_{t}\cdot{\gradt{\overline{q}}}}dxdy,\\
\mathcal{C}_{2}:&&\H(\mbox{curlt},\Gamma)\times{H^1(\Gamma)}\rightarrow{\mathbb{C}}\\
&&(\h_{t},q)\rightarrow\int_{\Gamma}{\d{\ep}_{t}^{-1}\h_{t}\cdot{\gradt{\overline{q}}}}dxdy.
\end{eqnarray*}
Using scalar Green's formula (\ref{scalargreen2}) and vector Green's formula (\ref{vectorgreen}), we get the variational forms of PDEs (\ref{3eqs2}) and (\ref{3eqsm2}), respectively.\\
\indent For the TE modes, find $k_{t}\in{\mathbb{R}}$, $0\neq \ee_t\in{\H_{0}(\mbox{curlt},\Gamma)}$, such that
\begin{subequations} \label{3eqsmeet1}
\begin{align}
\mathcal{A}_{1}(\ee_{t},\F)&= k_{t}^2 \mathcal{B}_{1}(\ee_{t},\F) ~~\forall~\F\in{\H_{0}(\text{curlt},\Gamma)},\label{formee1}\\
 \mathcal{C}_{1}(\ee_{t},q) &= 0 ~~\forall~q\in{H_{0}^1(\Gamma)}.\label{formee2}
\end{align}
\end{subequations}
\indent For the TM modes, find $k_{t}\in{\mathbb{R}}$, $0\neq \h_t\in{\H(\mbox{curlt},\Gamma)}$, such that
\begin{subequations} \label{3eqsmhht1}
\begin{align}
\mathcal{A}_{2}(\h_{t},\F)&= k_{t}^2 \mathcal{B}_{2}(\h_{t},\F) ~~\forall~\F\in{\H(\text{curlt},\Gamma)},\label{formhh1}\\
 \mathcal{C}_{2}(\h_{t},q) &= 0 ~~\forall~q\in{H^1(\Gamma)}.\label{formhh2}
\end{align}
\end{subequations}
It is clear that the variational forms (\ref{3eqsmeet1}) and (\ref{3eqsmhht1}) are equivalent to PDEs (\ref{3eqs2}) and (\ref{3eqsm2}), respectively. We must solve variational forms (\ref{formee1}) and (\ref{formhh1}) under the constraint of (\ref{formee2}) and (\ref{formhh2}), respectively. Based on Kikuchi's work \cite{Kikuchi}, we try to change the variation forms (\ref{3eqsmeet1}) and (\ref{3eqsmhht1}) to the corresponding mixed variational forms by means of a Lagrangian multiplier. \\
\indent Let us now introduce the mixed variational form associated with the TE modes:
Seek $k_{t}\in{\mathbb{R}}$,  $0\neq\ee_{t}\in{\H_{0}(\mbox{curlt},\Gamma)}$, $p_{1}\in{H_{0}^1(\Gamma)}$ such that
\begin{subequations} \label{3eqsmeetmix1}
\begin{align}
\mathcal{A}_{1}(\ee_{t},\F)+\overline{\mathcal{C}}_{1}(\F,p_{1})&= k_{t}^2 \mathcal{B}_{1}(\ee_{t},\F)~\forall~\F\in{\H_{0}(\text{curlt},\Gamma)},\label{formeemix1}\\
 \mathcal{C}_{1}(\ee_{t},q) &= 0 ~~\forall~q\in{H_{0}^1(\Gamma)}.\label{formeemix2}
\end{align}
\end{subequations}
\indent Let us now introduce the mixed variational form associated with the TM modes:
Seek $k_{t}\in{\mathbb{R}}$,  $0\neq\h_{t}\in{\H(\mbox{curlt},\Gamma)}$, $p_{2}\in{H^1(\Gamma)}$ such that
\begin{subequations} \label{3eqsmhhtmix1}
\begin{align}
\mathcal{A}_{2}(\h_{t},\F)+\overline{\mathcal{C}}_{2}(\F,p_{2})&= k_{t}^2 \mathcal{B}_{2}(\h_{t},\F)~\forall~\F\in{\H(\text{curlt},\Gamma)},\label{formhhmix1}\\
 \mathcal{C}_{2}(\h_{t},q) &= 0 ~~\forall~q\in{H^1(\Gamma)},\label{formhhmix2}
\end{align}
\end{subequations}
where $p_{k}~(k=1,2)$ are the Lagrangian multipliers and $\overline{\mathcal{C}}_{k}(\F,p_{k})~(k=1,2)$ stand for the complex conjugation of the continuous sesquilinear form $\mathcal{C}_{k}(\F,p_{k})~(k=1,2)$. The weak forms (\ref{3eqsmeetmix1}) and (\ref{3eqsmhhtmix1}) are the saddle point problems in finite element analysis. \\
\indent We have already proved the equivalence between the variational form (\ref{3eqsmeet1}) and mixed variational form (\ref{3eqsmeetmix1}) in \cite{jiang2014}, because  $p_{1}=0\mbox{~in~}\Gamma$ is proved in \cite{jiang2014}. Now we prove the equivalence between the variational form (\ref{3eqsmhht1}) and
mixed variational form (\ref{3eqsmhhtmix1}). Obviously, any eigenpair of (\ref{3eqsmhht1}) with $p_{2}=C$ satisfies (\ref{3eqsmhhtmix1}), where $C$ is an arbitrary constant.
Conversely, by taking $\F=\gradt{p_{2}}$ in (\ref{formhhmix1}) and $q=p_{2}$ in (\ref{formhhmix2}), we get
$\overline{\mathcal{C}}_{2}(\gradt{p_2},p_2)=k_{t}^2 \mathcal{B}_{2}(\h_{t},\gradt{p_2})=k_{t}^2 \mathcal{C}_{2}(\h_{t},\gradt{p_2})=0$.
Since $\d{\ep}_{t}^{-1}$ has positive real eigenvalues, we can deduce that $\gradt{p_2}=0$, therefore $p_2=C$ and $\overline{\mathcal{C}}_{2}(\F,p_{2})=0$, which shows that any eigenpair of the mixed variational form (\ref{3eqsmhhtmix1}) satisfies (\ref{3eqsmhht1}) as well. Numerical experiments are also shown that $p_2=C\mbox{~in~}\Gamma$. Next we utilize mixed FEM to discretize the mixed variational problems (\ref{3eqsmeetmix1}) and (\ref{3eqsmhhtmix1}), respectively.\\
\indent Let $\mathcal{T}_{h}$ be a regular triangular partition \cite{ciar} of $\Gamma$ with mesh parameter $h$.
About the definition of CT/LN edge element space, please see \cite{jiang2014}. Here we only consider the discretization of mixed variation form (\ref{3eqsmhhtmix1}). For
the discretization of mixed variation form (\ref{3eqsmeetmix1}),
for details, please see \cite{jiang2014}.
Let $\W_{h}^{2}$ be CT/LN edge element space and $\W_{h}^{2}=\mbox{span}\{\N_{1}, \N_{2},...,\N_{s}\}$, where $\N_{k}$ is the $k$-th global basis function in $\W_{h}^{2}$ and the integer $s$ is the number of total edges (including the edges on the boundary) in triangular mesh $\mathcal{T}_{h}$.
From \cite{Petermonk,Hiptmair}, we know that $\W_{h}^{2}$ is a FEM subspace of $\H(\text{curlt},\Gamma)$. Let $S_{h}^{2}$ be the standard linear element space and $S_{h}^{2}=\mbox{span}\{\phi_{1}, \phi_{2},...,\phi_{n}\}$, where $\phi_{k}$ is the $k$-th  global basis function in $S_{h}^{2}$ and the integer $n$ is the number of total nodes (including the nodes on the boundary) in triangular mesh $\mathcal{T}_{h}$. From \cite{ciar}, we know that $S_{h}^{2}$ is a FEM subspace of $H^{1}(\Gamma)$. \\
\indent Restricting the mixed variational form (\ref{3eqsmhhtmix1}) on $\W_{h}^{2}\times S_{h}^{2}$, we obtain the discrete mixed variational form:
Seek $k_{t,h}\in{\mathbb{R}}$,  $0\neq\h_{t}^{h}\in{\W_{h}^{2}}$, $p_{2,h}\in{S_{h}^{2}}$ such that
\begin{subequations} \label{3eqsmhhtmixss1}
\begin{align}
\mathcal{A}_{2}(\h_{t}^{h},\F)+\overline{\mathcal{C}}_{2}(\F,p_{2,h})&= k_{t,h}^2 \mathcal{B}_{2}(\h_{t}^{h},\F)~\forall~\F\in{\W_{h}^{2}},\label{formhhmixss1}\\
 \mathcal{C}_{2}(\h_{t}^{h},q) &= 0 ~~\forall~q\in{S_{h}^{2}}.\label{formhhmixss2}
\end{align}
\end{subequations}
The discrete weak form (\ref{3eqsmhhtmixss1}) is a conforming mixed FEM discretization of mixed variational form (\ref{3eqsmhhtmix1}).
Since $\h_{t}^{h}\in{W_{h}^{2}}$ and $p_{2,h}\in{S_{h}^2}$, we write
$$\h_{t}^{h}=\sum_{k=1}^{s}\xi_{k}\N_{k},\quad p_{2,h}=\sum_{k=1}^{n}\zeta_{k}\phi_{k}.$$
Finally, we get the following generalized matrix eigenvalue problem:
\begin{equation}\label{eigm1}
    \left[
      \begin{array}{cc}
        A_{2} & C_{2} \\
        C_{2}^{\dagger} & O \\
      \end{array}
    \right]\left[
             \begin{array}{c}
               \xi \\
               \zeta \\
             \end{array}
           \right]
    =k_{t,h}^2\left[
      \begin{array}{cc}
        B_{2} & O \\
        O & O \\
      \end{array}
    \right]\left[
             \begin{array}{c}
               \xi \\
               \zeta \\
             \end{array}
           \right],
\end{equation}
where
\begin{gather*}
   \xi=[\xi_{1},\xi_{2},\cdots,\xi_{s}]^{T}, \quad\zeta=[\zeta_{1},\zeta_{2},\cdots,\zeta_{n}]^{T},\\
   A_{2}=(a_{ij}^{(2)})\in{\mathbb{C}^{s\times s}},~C_{2}=(c_{ij}^{(2)})\in{\mathbb{C}^{s\times n}},\\
   B_{2}=(b_{ij}^{(2)})\in{\mathbb{C}^{n\times n}},~a_{ij}^{(2)}=\mathcal{A}_{2}(\N_{j},\N_{i}),\\
   c_{ij}^{(2)}=\overline{\mathcal{C}}_2(\N_{i},\phi_{j}),~b_{ij}^{(2)}=\mathcal{B}_{2}(\N_{j},\N_{i}).
\end{gather*}
After solving the generalized algebraic eigenvalue problem (\ref{eigm1}), we can get the distribution of the transverse magnetic field $\h_{t}^{h}$ in $\Gamma$ using an interpolation technique.\\
\indent For the implementation of the mixed variational form (\ref{3eqsmhhtmixss1}),
we do not need to deal with the boundary conditions,
because two boundary conditions in PDE (\ref{3eqsm2}) are all natural boundary conditions in FEM. However, about the numerical implementation of the mixed variational form (\ref{3eqsmeetmix1}), the numerical treatment method is almost same as the mixed variational form (\ref{3eqsmhhtmix1}), but at last we need to enforce boundary condition, because for the TE modes, $\^n\times\ee_{t}=0\mbox{~on~}\Gamma$ is an essential boundary condition in FEM.
\section{Numerical Experiments}
\begin{figure}[ht]
  \centering
  {\subfigure[]{
    \label{aaa}
   \includegraphics[width=0.45\columnwidth,draft=false]{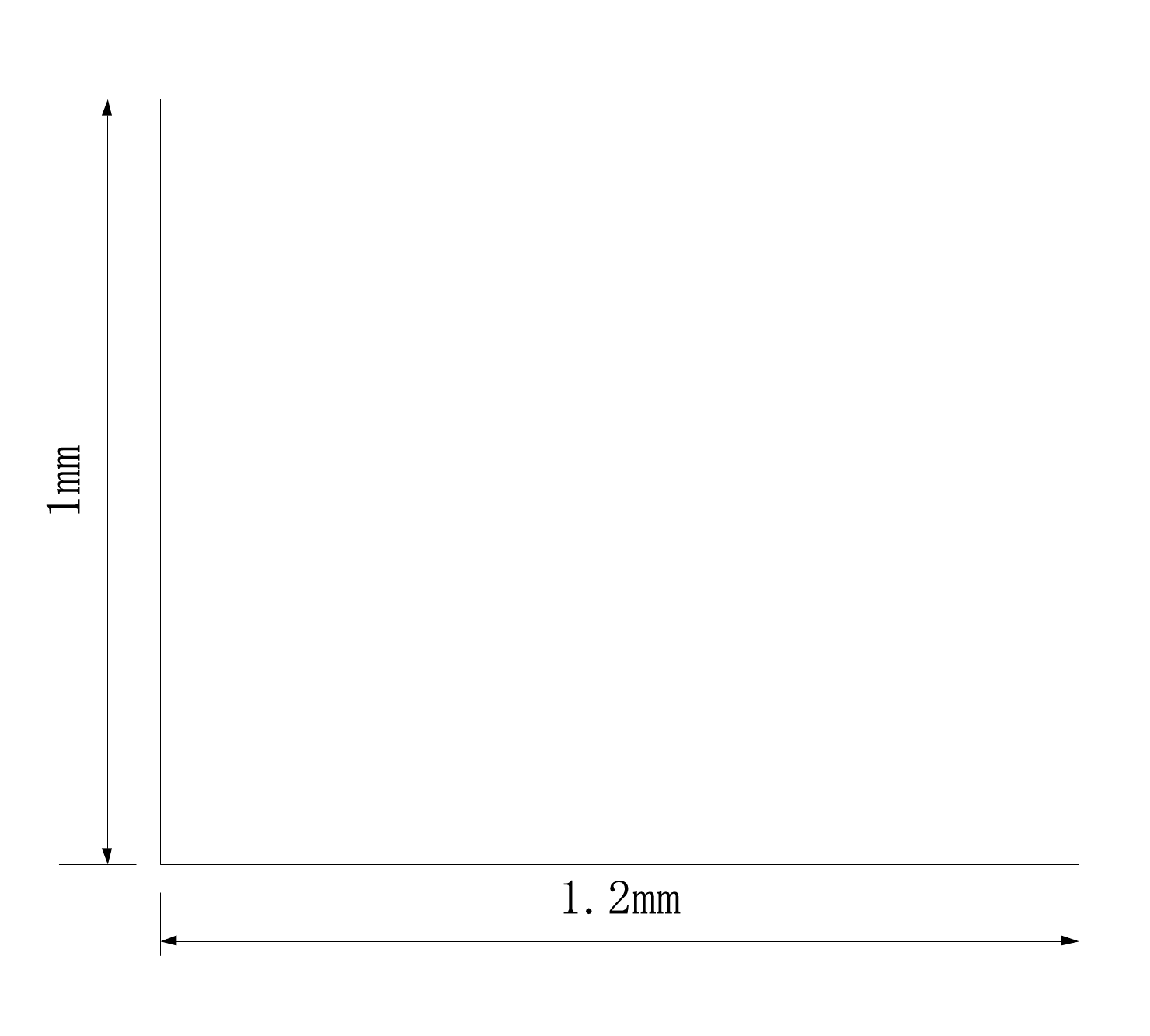}}
     \subfigure[]{
    \label{bbbph}
   \includegraphics[width=0.4\columnwidth,draft=false]{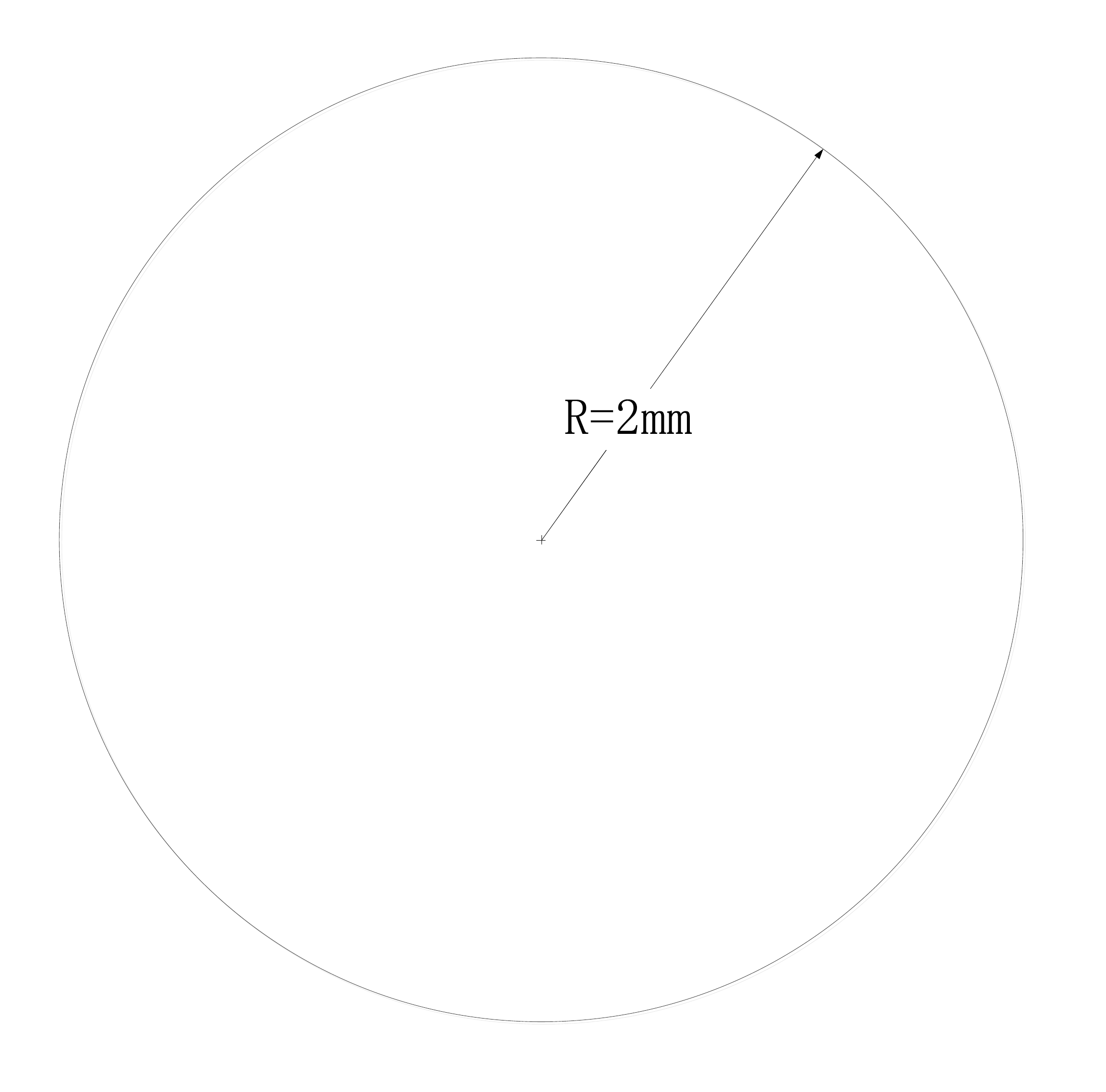}}\\
   \subfigure[]{
    \label{ccc}
   \includegraphics[width=0.4\columnwidth,draft=false]{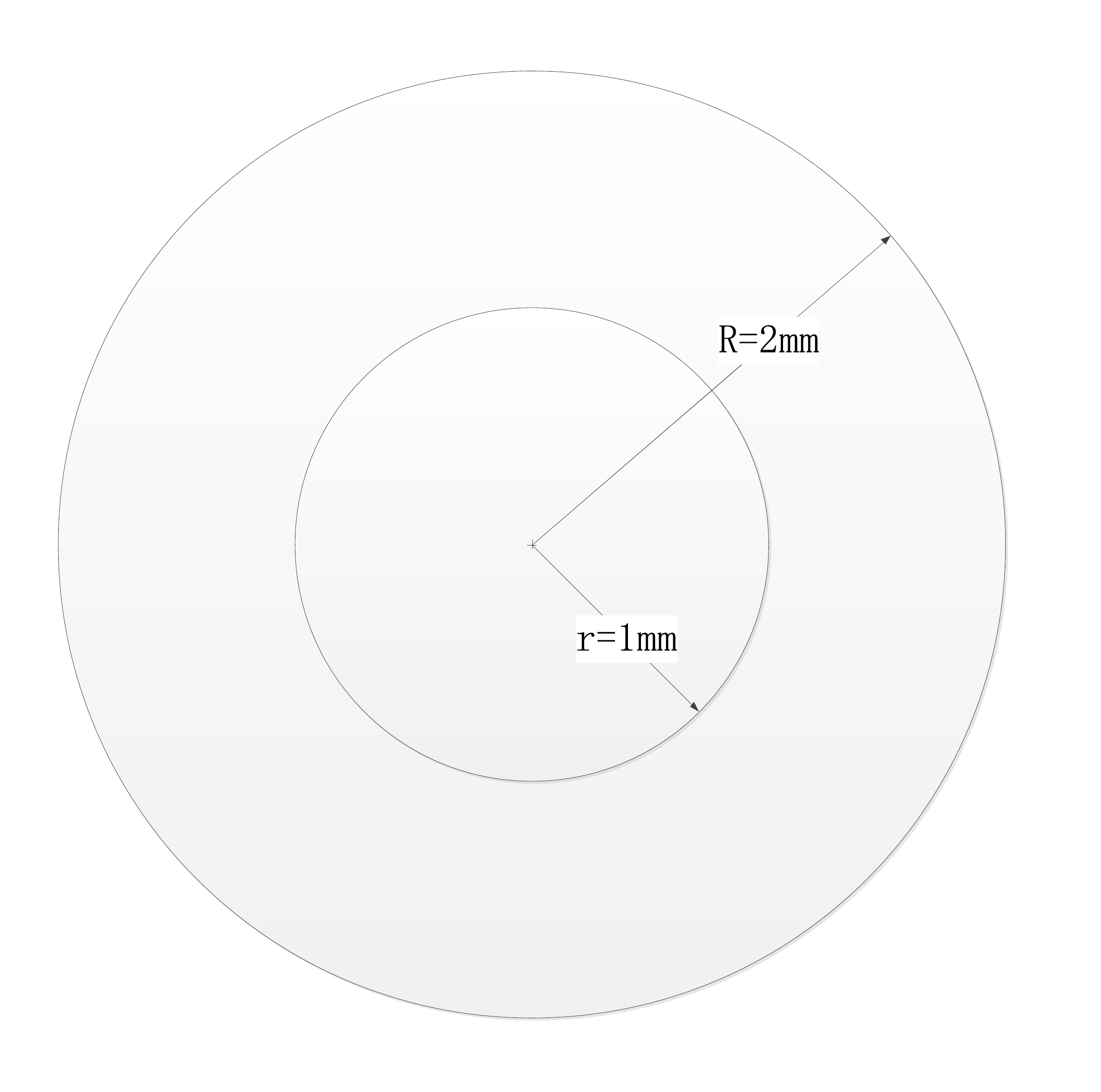}}
   \subfigure[]{
    \label{ddd}
   \includegraphics[width=0.48\columnwidth,draft=false]{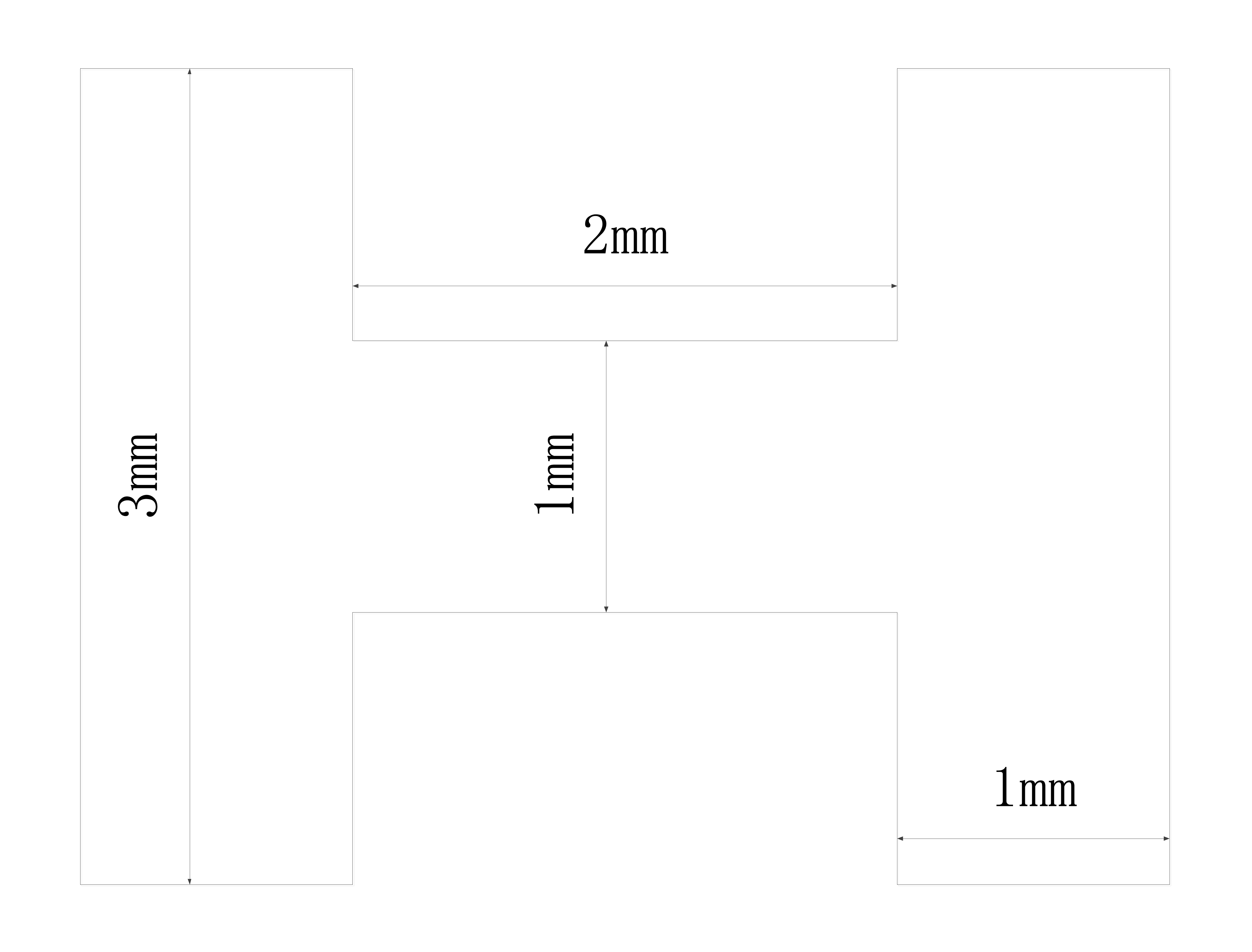}}}
 \caption{(a) Dimension of Rectangular waveguide. (b) Dimension of cylindrical waveguide. (c) Dimension of coaxial waveguide. (d) Dimension of double-ridge waveguide.}
\label{figconv}
\end{figure}

\begin{figure}[ht]
  \centering
  {\subfigure[]{
    \label{terecreal1}
   \includegraphics[width=0.48\columnwidth,draft=false]{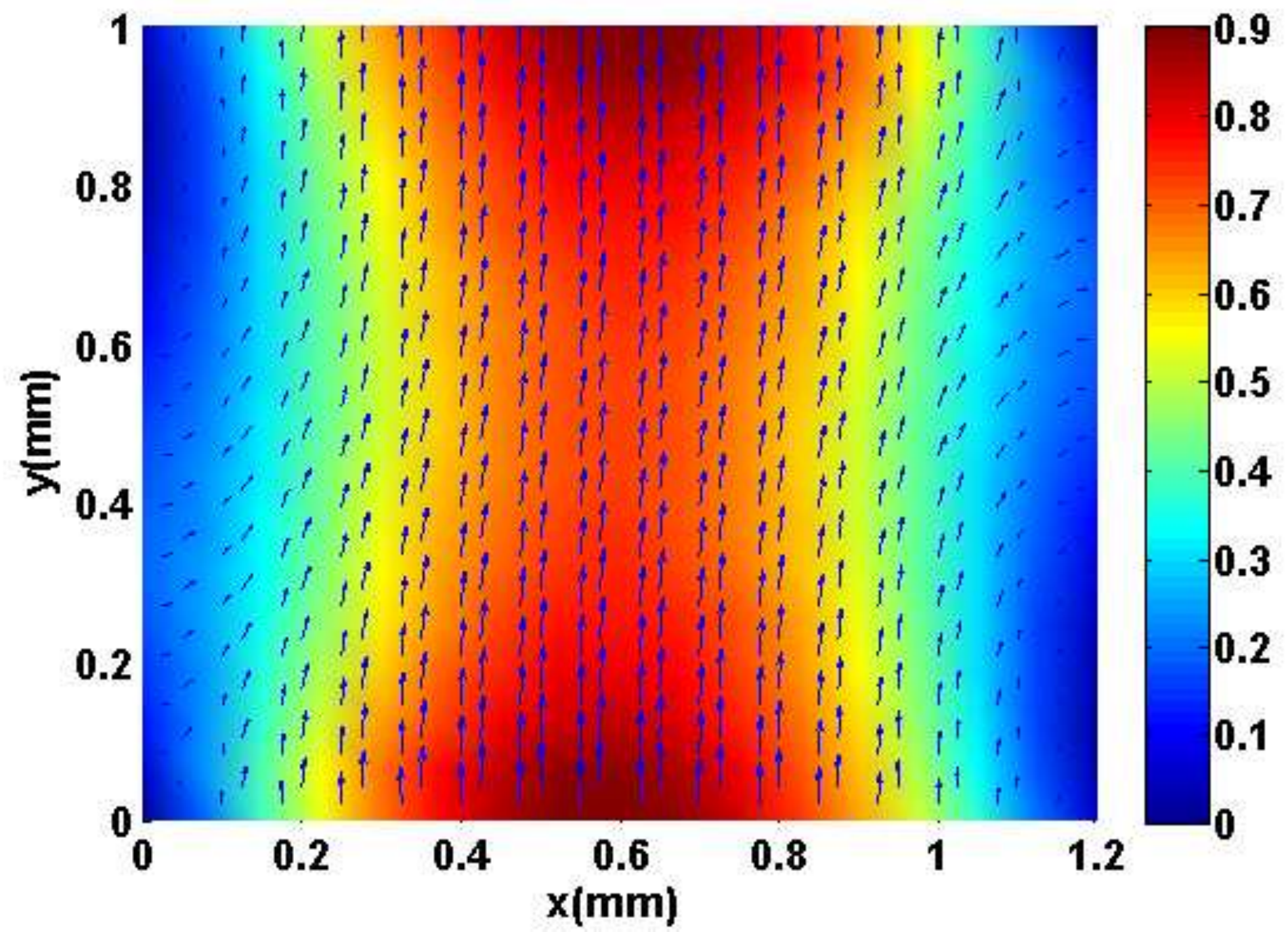}}
     \subfigure[]{
    \label{terecreal2}
   \includegraphics[width=0.48\columnwidth,draft=false]{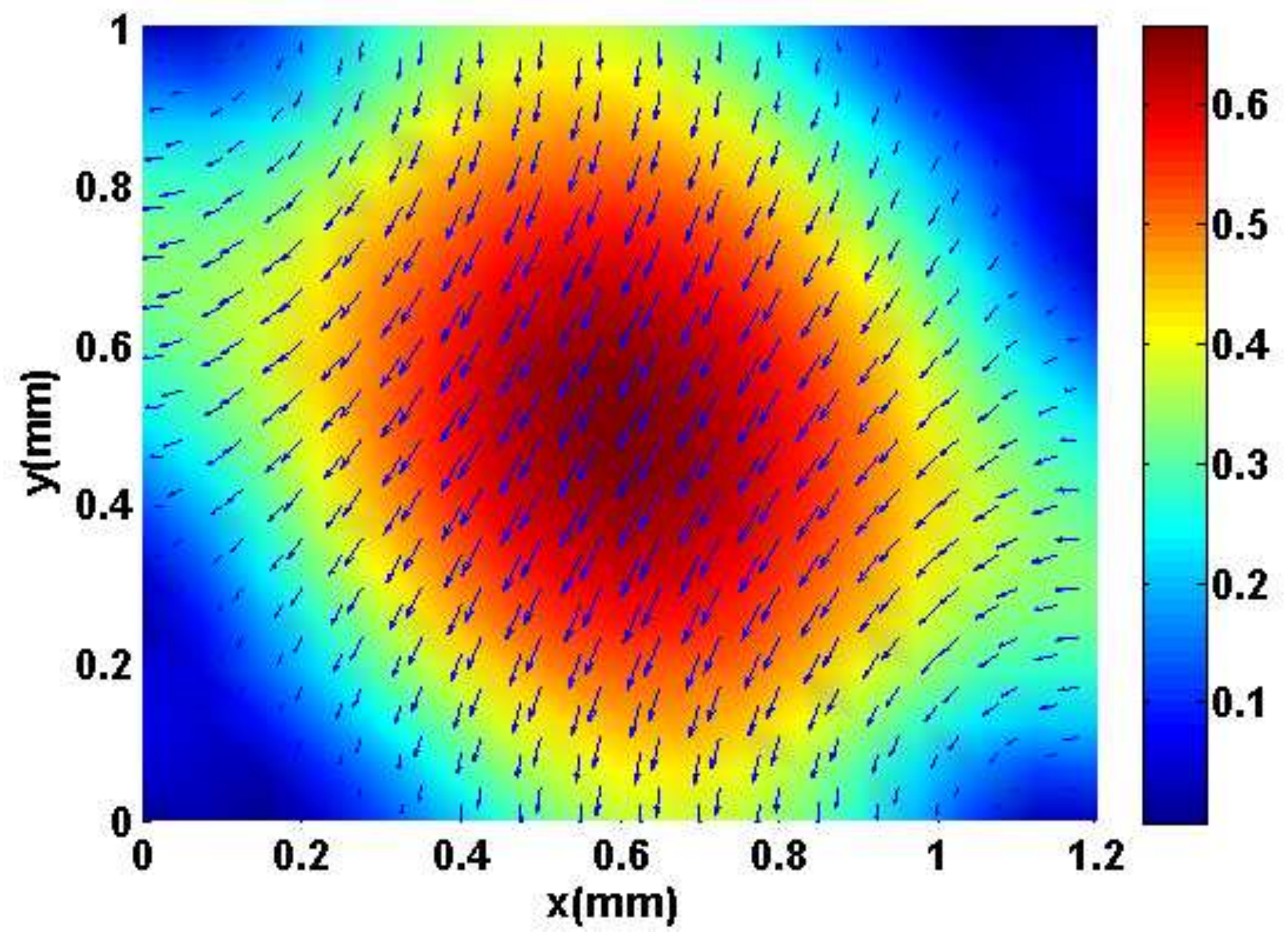}}\\
   \subfigure[]{
    \label{terecreal3}
   \includegraphics[width=0.48\columnwidth,draft=false]{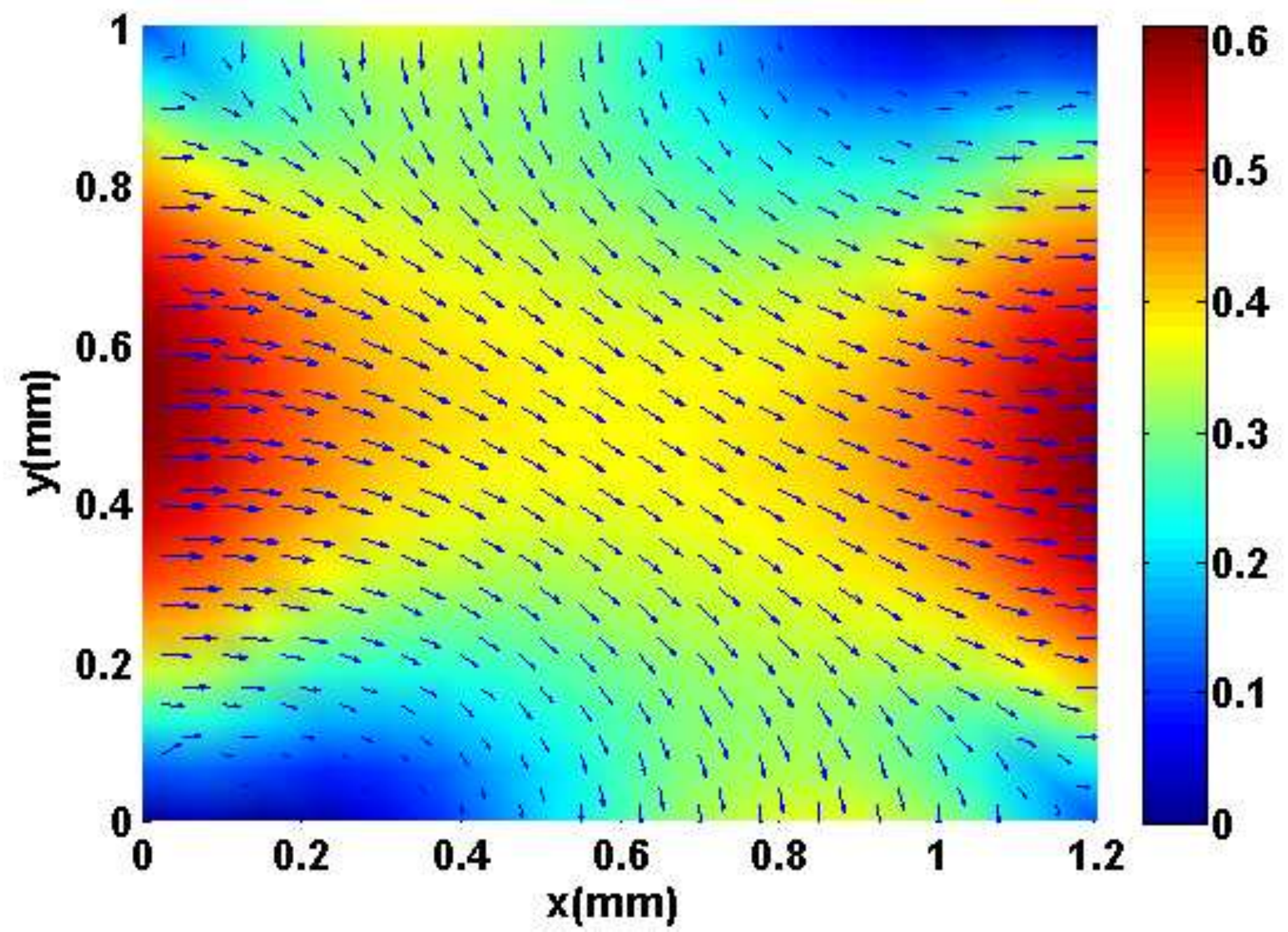}}
   \subfigure[]{
    \label{terecreal4}
   \includegraphics[width=0.48\columnwidth,draft=false]{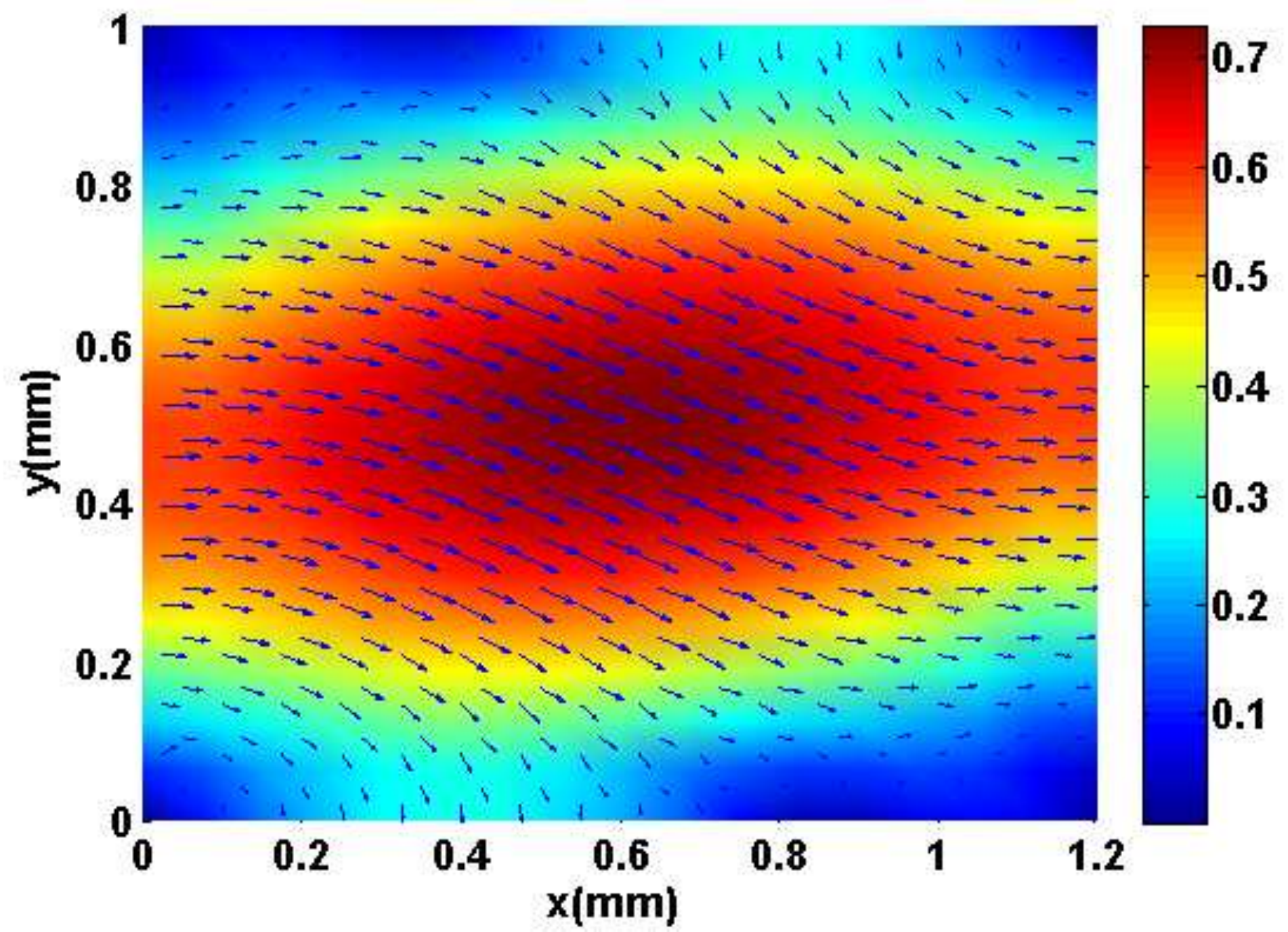}}\\
      \subfigure[]{
    \label{terecreal5}
   \includegraphics[width=0.48\columnwidth,draft=false]{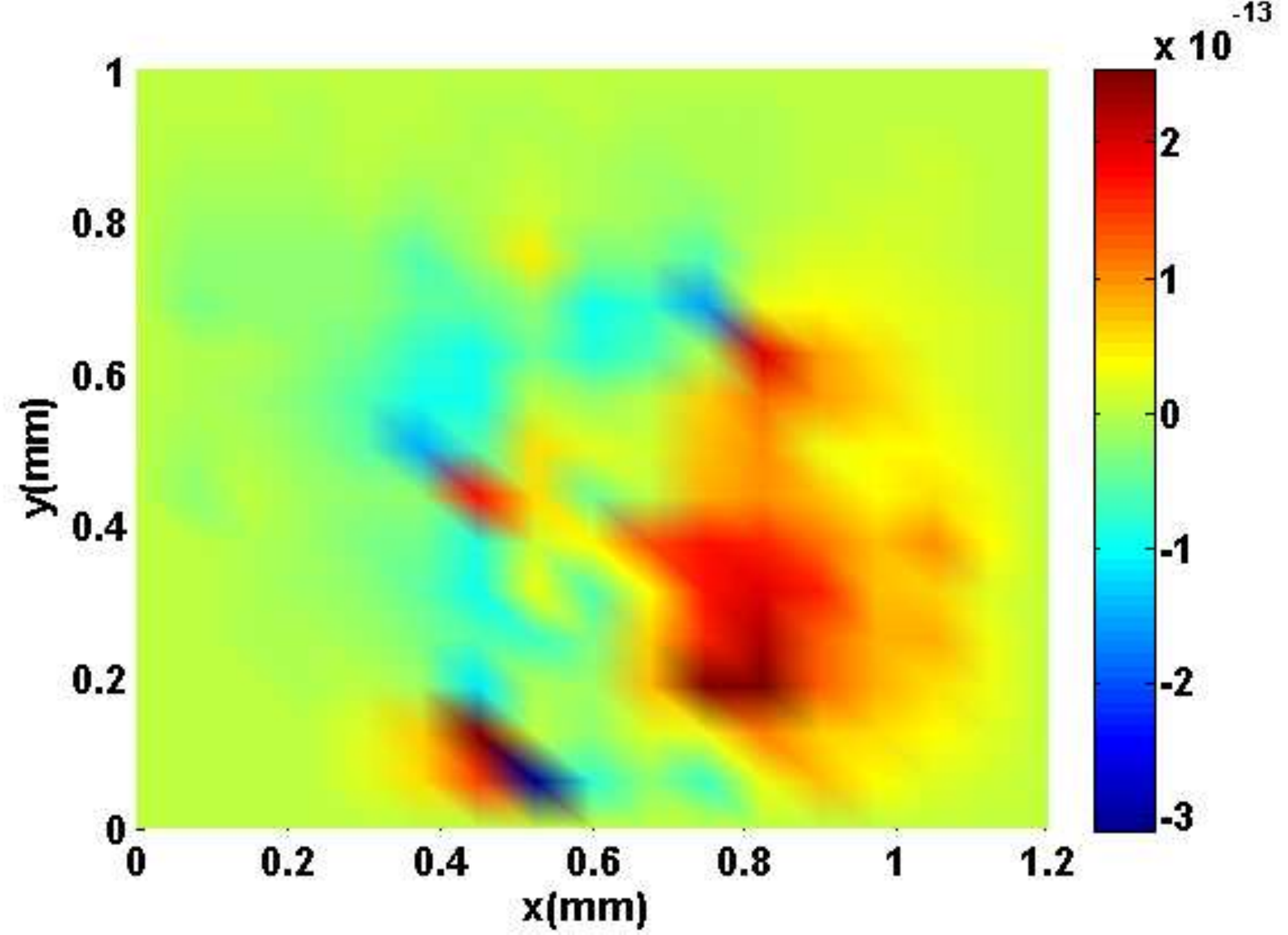}}
   \subfigure[]{
    \label{terecreal6}
   \includegraphics[width=0.48\columnwidth,draft=false]{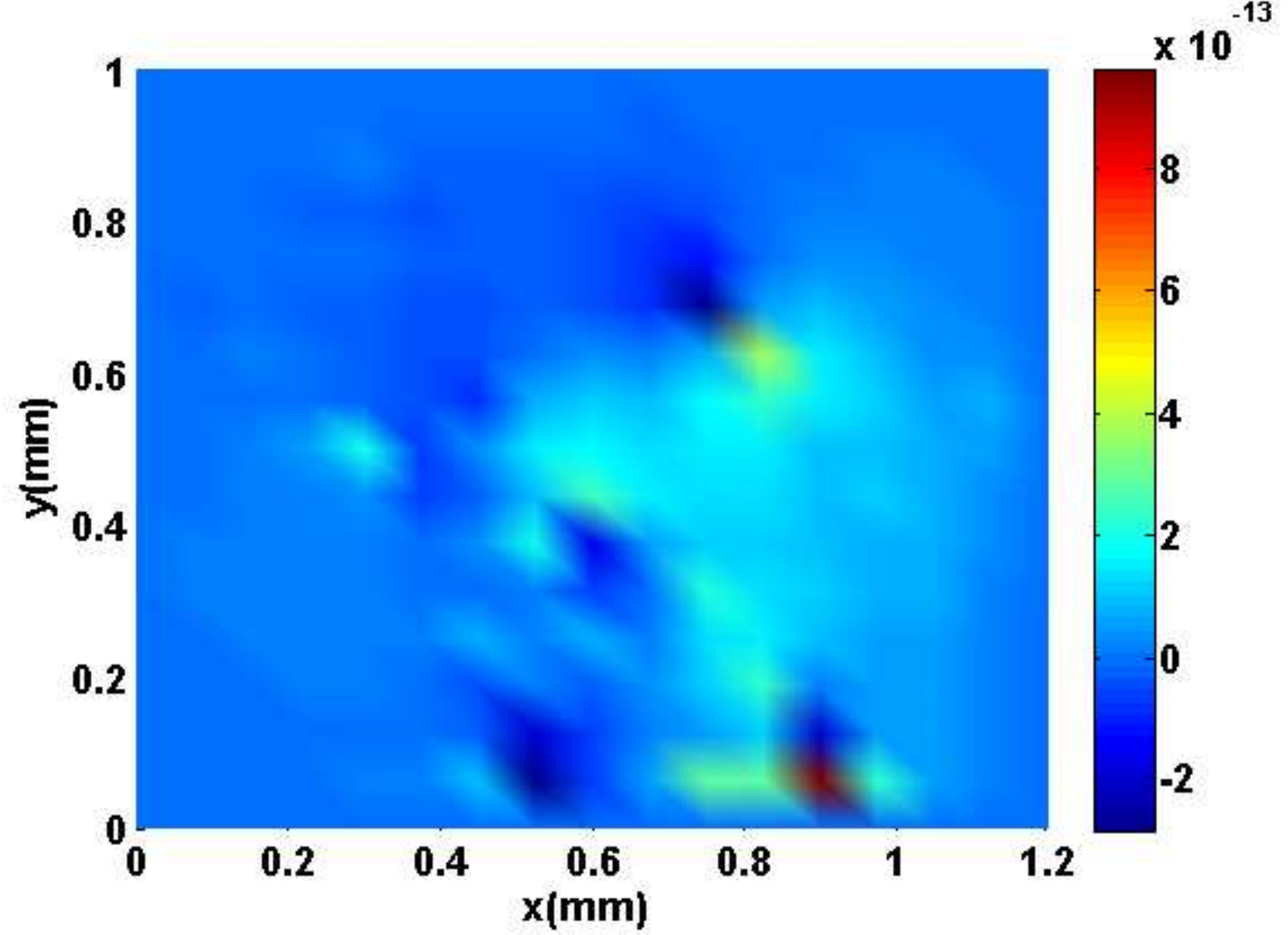}}\\}
 \caption{(a) Magnitude and field distribution of $\mbox{Re}(\ee_{t,1}^{h})$.
          (b) Magnitude and field distribution of $\mbox{Re}(\ee_{t,2}^{h})$.
          (c) Magnitude and field distribution of $\mbox{Im}(\ee_{t,1}^{h})$.
          (d) Magnitude and field distribution of $\mbox{Im}(\ee_{t,2}^{h})$.
          (e) Magnitude distribution of $p_{1,h}^{(1)}$.
          (f) Magnitude distribution of $p_{1,h}^{(2)}$.
          These figures are obtained from mixed FEM on the second mesh ($h\approx0.0976$) for the first two independent TE modes in rectangular waveguide.}
\label{terecwaveguide}
\end{figure}

\begin{figure}[ht]
  \centering
  {\subfigure[]{
    \label{tmrecreal1}
   \includegraphics[width=0.48\columnwidth,draft=false]{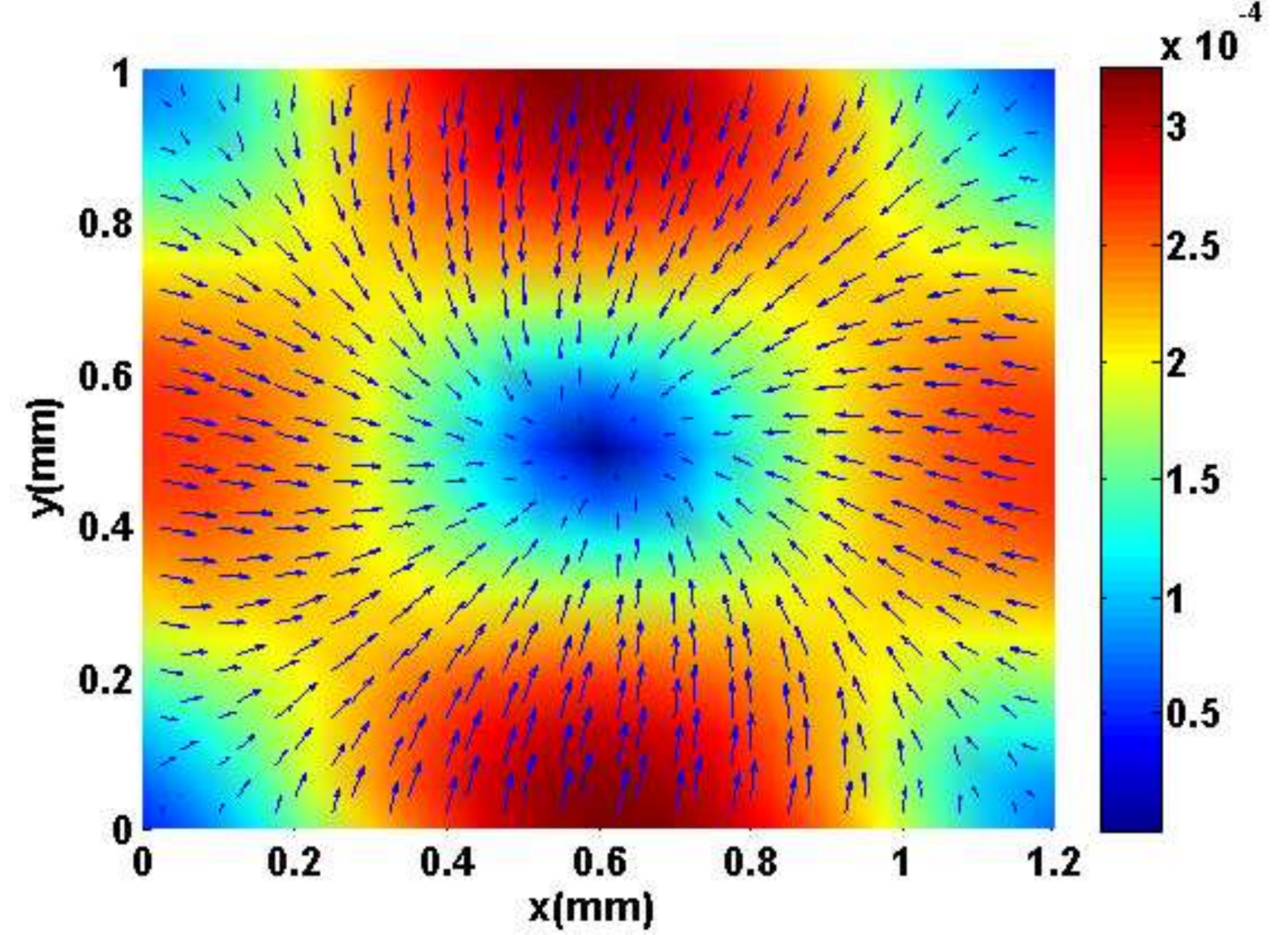}}
     \subfigure[]{
    \label{tmrecreal2}
   \includegraphics[width=0.48\columnwidth,draft=false]{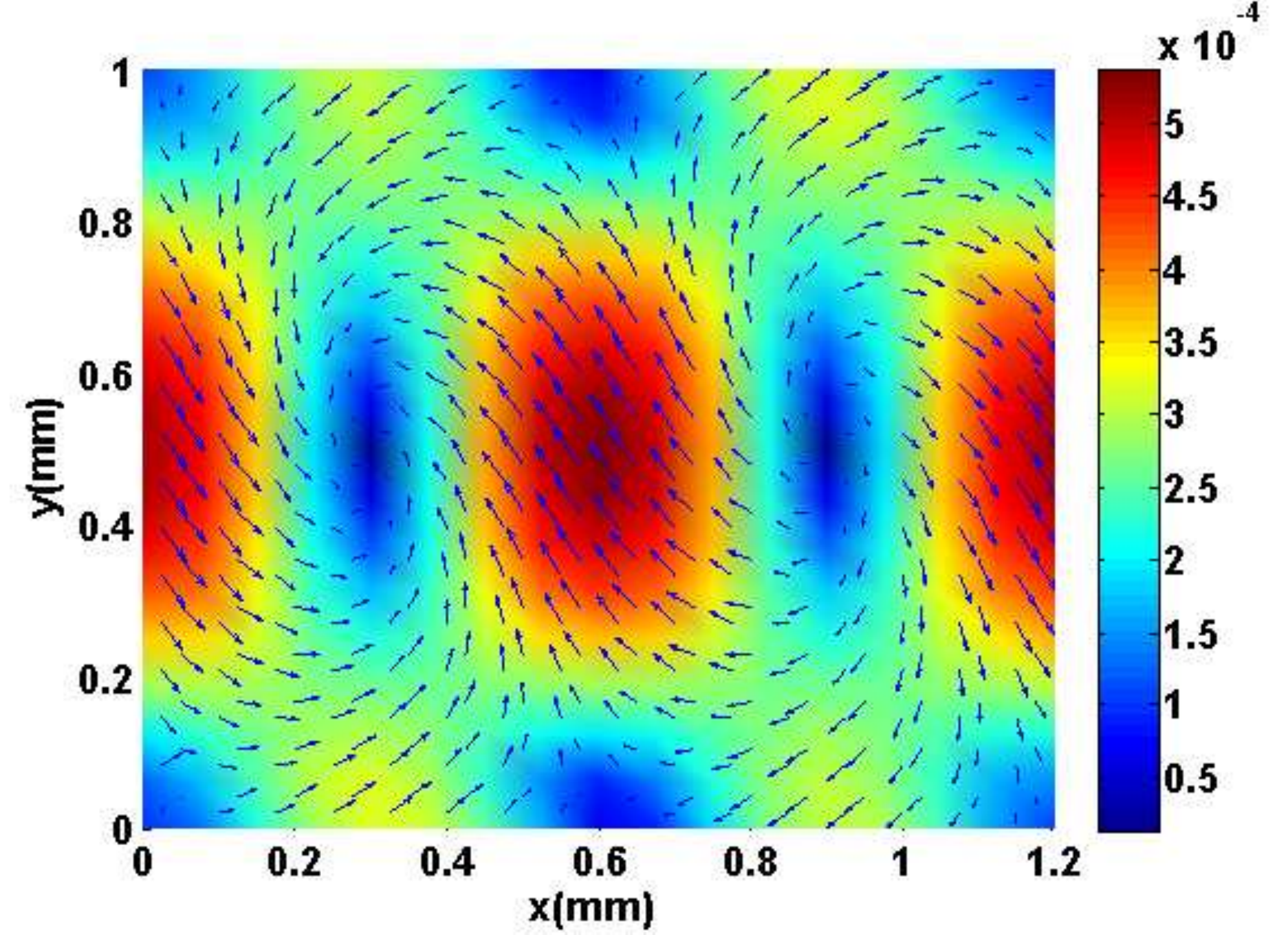}}\\
   \subfigure[]{
    \label{tmrecreal3}
   \includegraphics[width=0.48\columnwidth,draft=false]{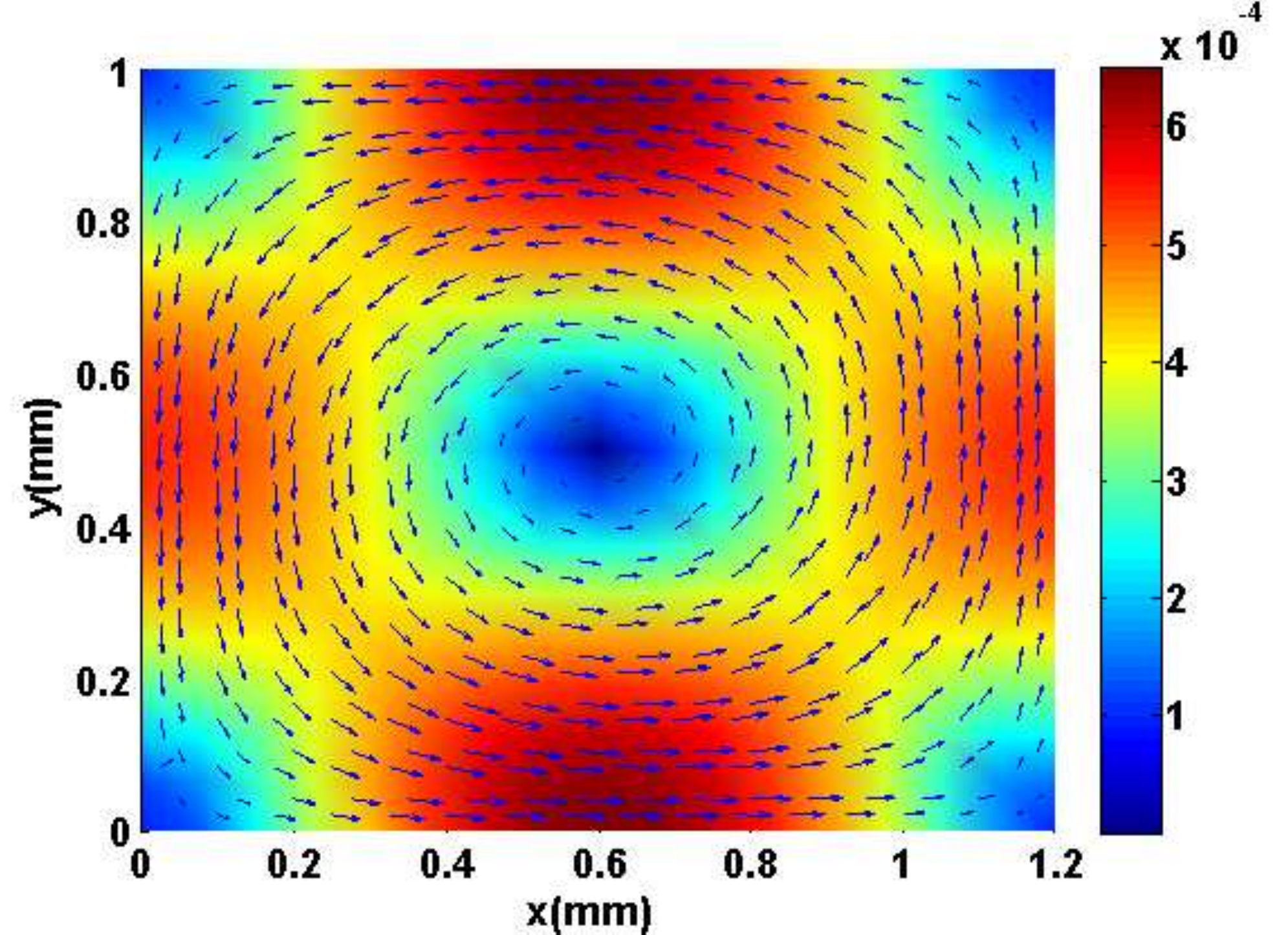}}
   \subfigure[]{
    \label{tmrecreal4}
   \includegraphics[width=0.48\columnwidth,draft=false]{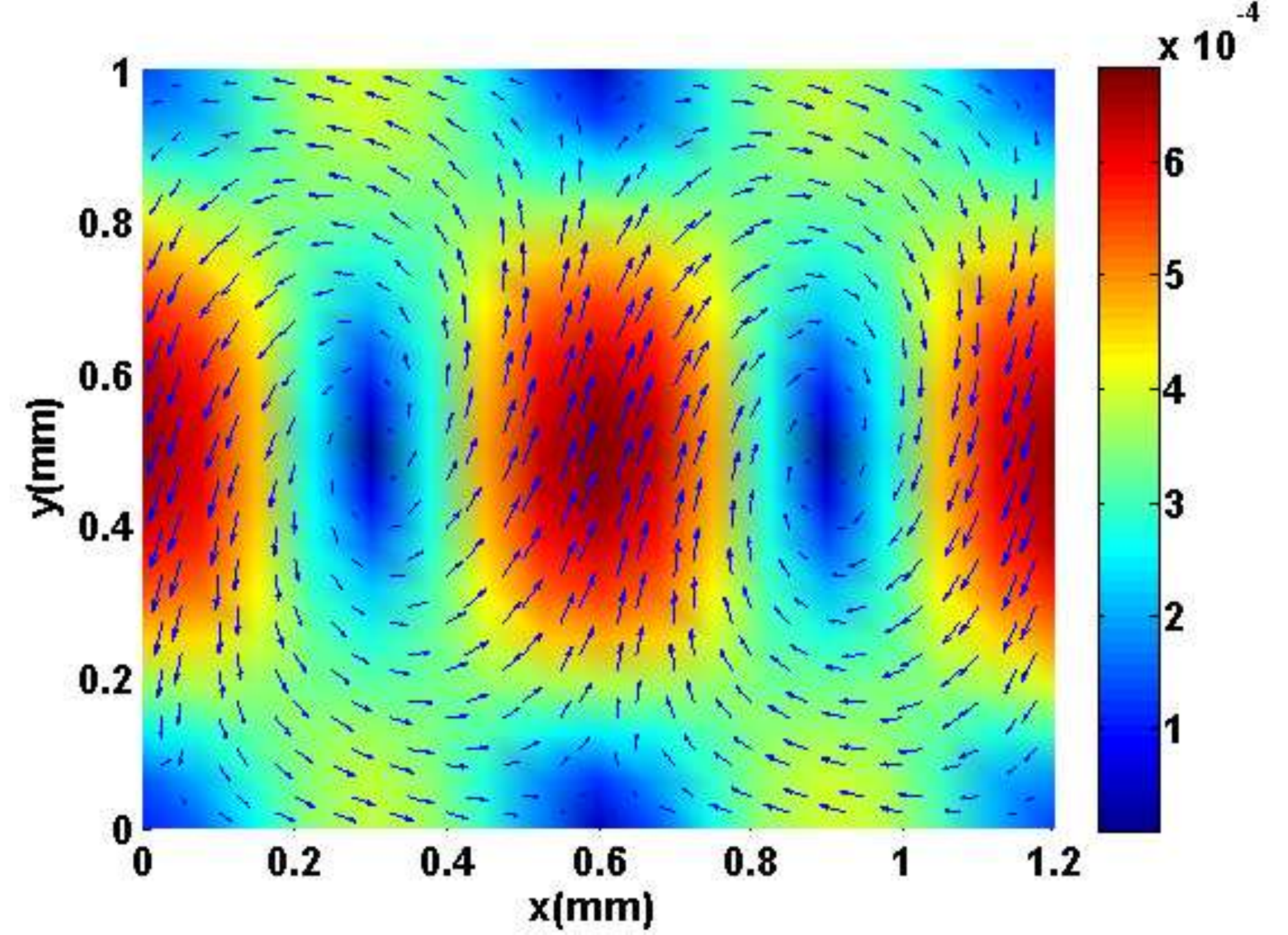}}\\
      \subfigure[]{
    \label{tmrecreal5}
   \includegraphics[width=0.48\columnwidth,draft=false]{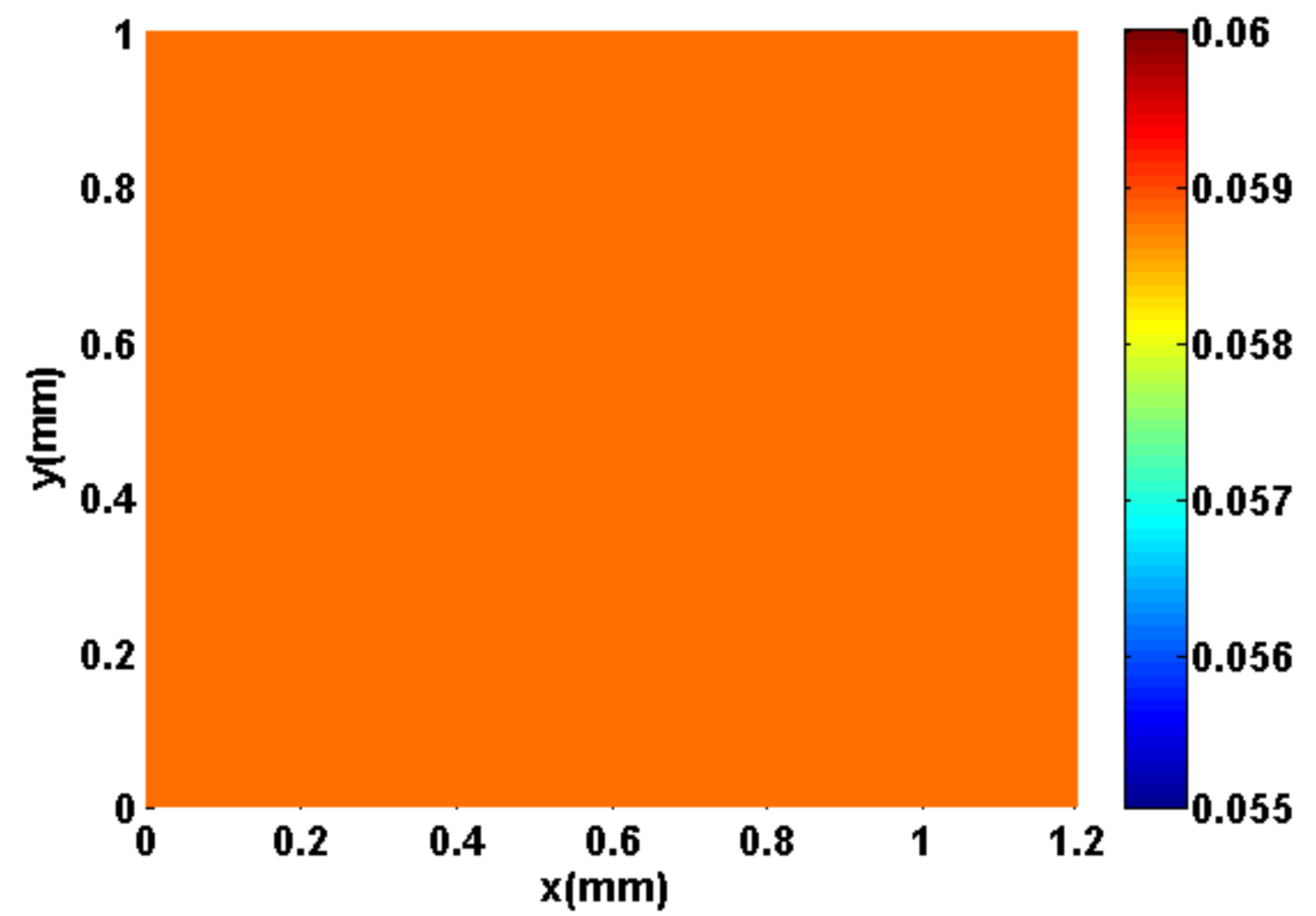}}
   \subfigure[]{
    \label{tmrecreal6}
   \includegraphics[width=0.48\columnwidth,draft=false]{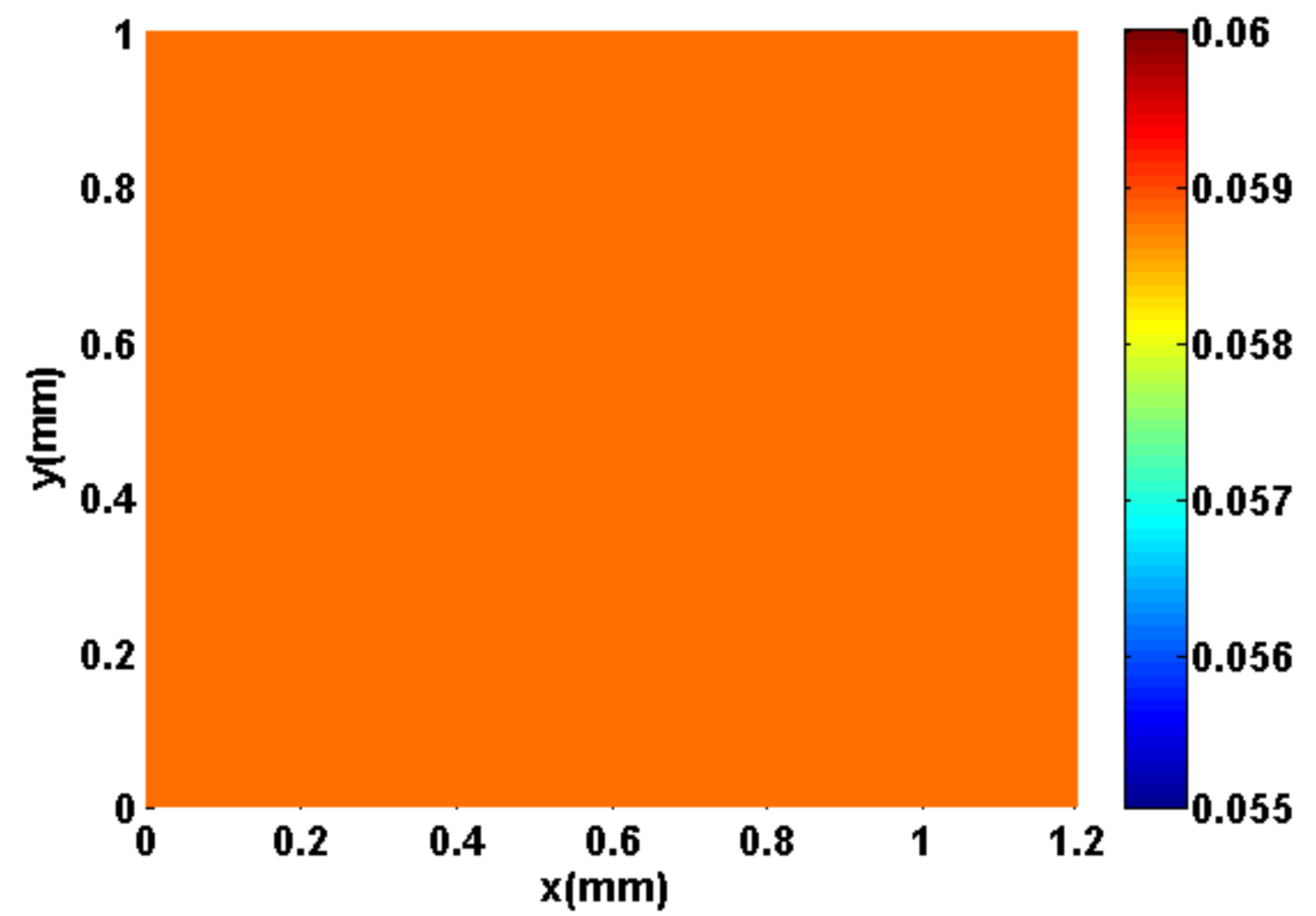}}\\}
 \caption{(a) Magnitude and field distribution of $\mbox{Re}(\h_{t,1}^{h})$.
          (b) Magnitude and field distribution of $\mbox{Re}(\h_{t,2}^{h})$.
          (c) Magnitude and field distribution of $\mbox{Im}(\h_{t,1}^{h})$.
          (d) Magnitude and field distribution of $\mbox{Im}(\h_{t,2}^{h})$.
          (e) Magnitude distribution of $p_{2,h}^{(1)}$.
          (f) Magnitude distribution of $p_{2,h}^{(2)}$.
          These figures are obtained from mixed FEM on the second mesh ($h\approx0.0976$) for the first two independent TM modes in rectangular waveguide.}
\label{tmrecwaveguide}
\end{figure}

\begin{figure}[ht]
  \centering
  {\subfigure[]{
    \label{tecoareal1}
   \includegraphics[width=0.48\columnwidth,draft=false]{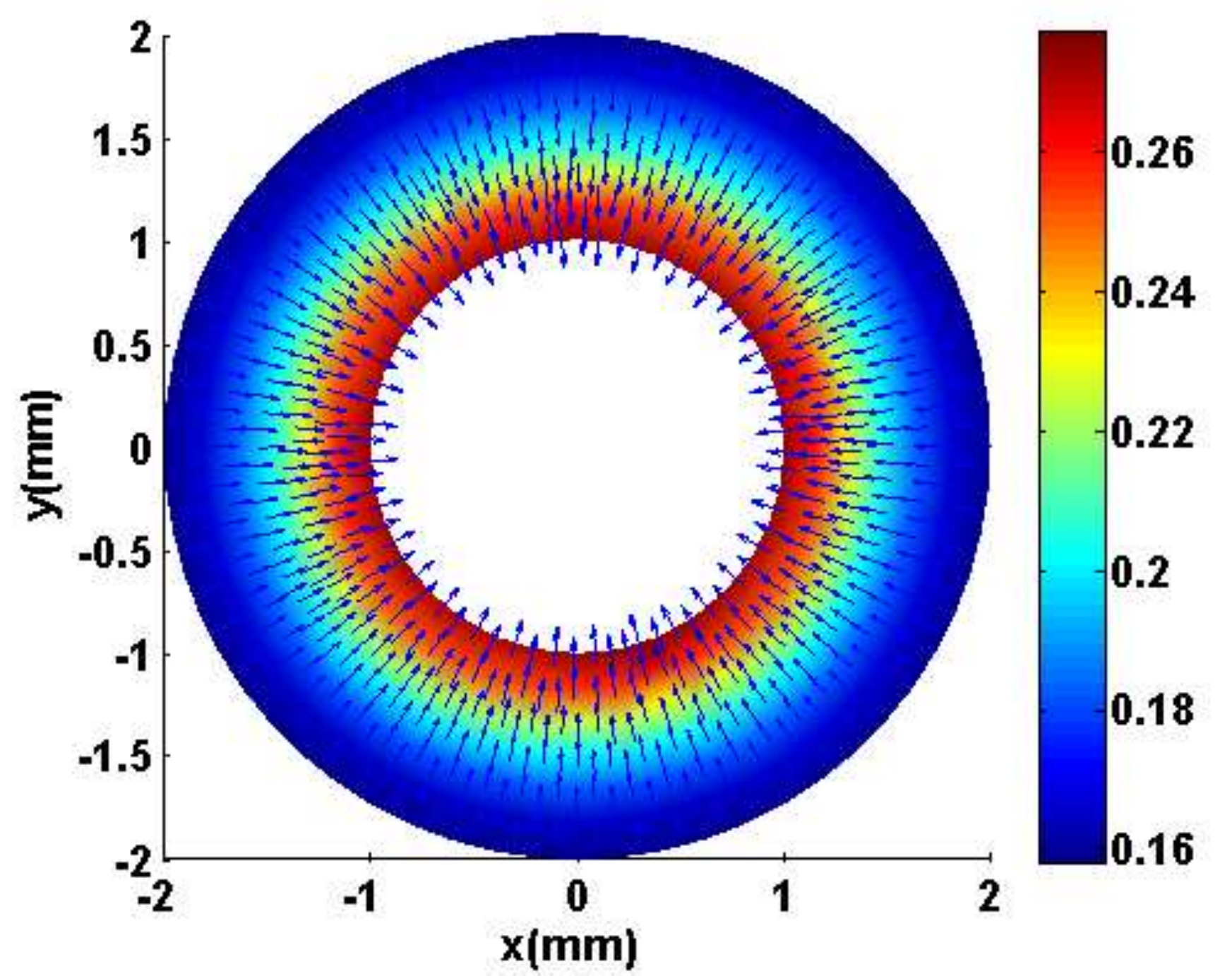}}
     \subfigure[]{
    \label{tecoareal2}
   \includegraphics[width=0.48\columnwidth,draft=false]{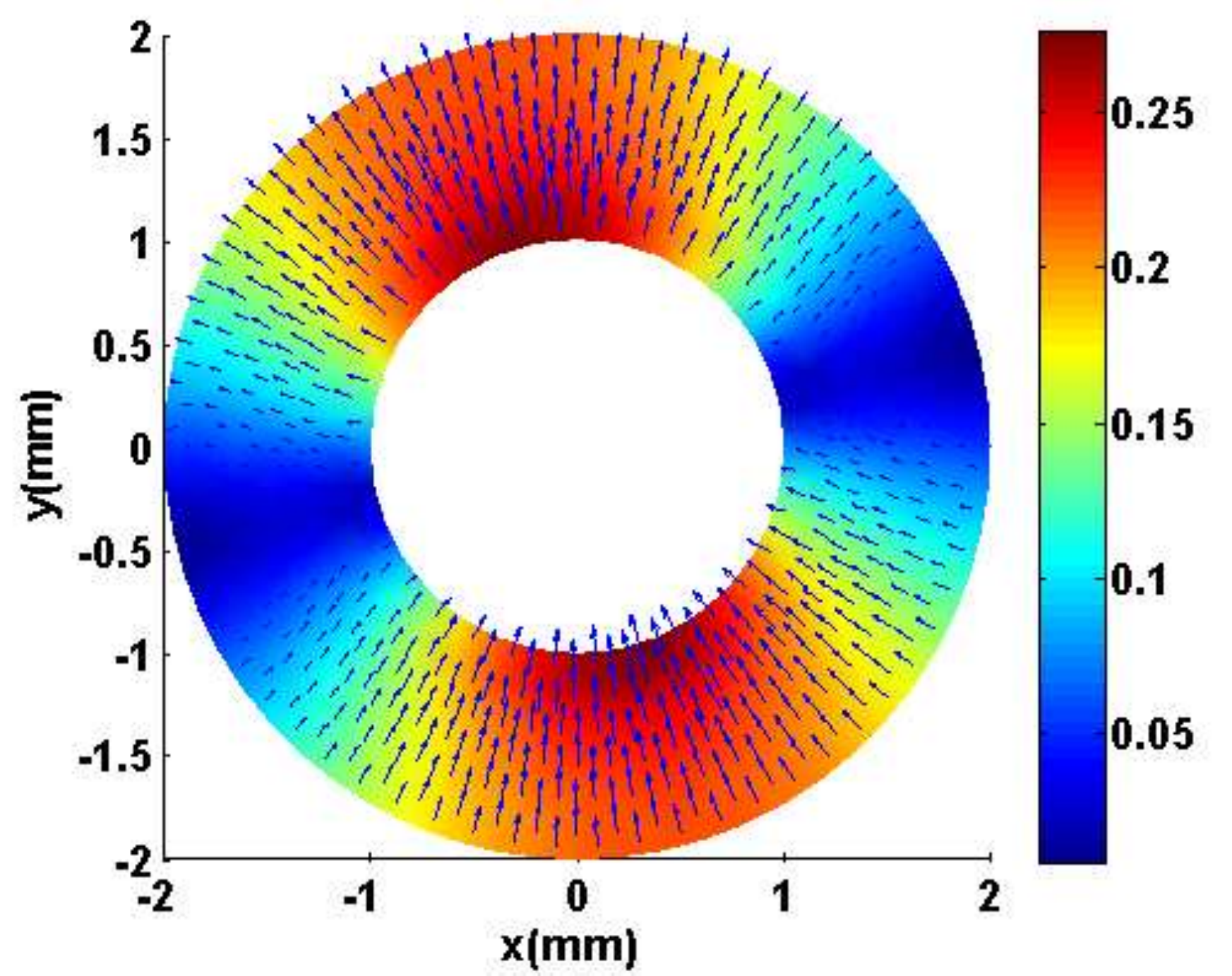}}\\
   \subfigure[]{
    \label{tecoareal3}
   \includegraphics[width=0.48\columnwidth,draft=false]{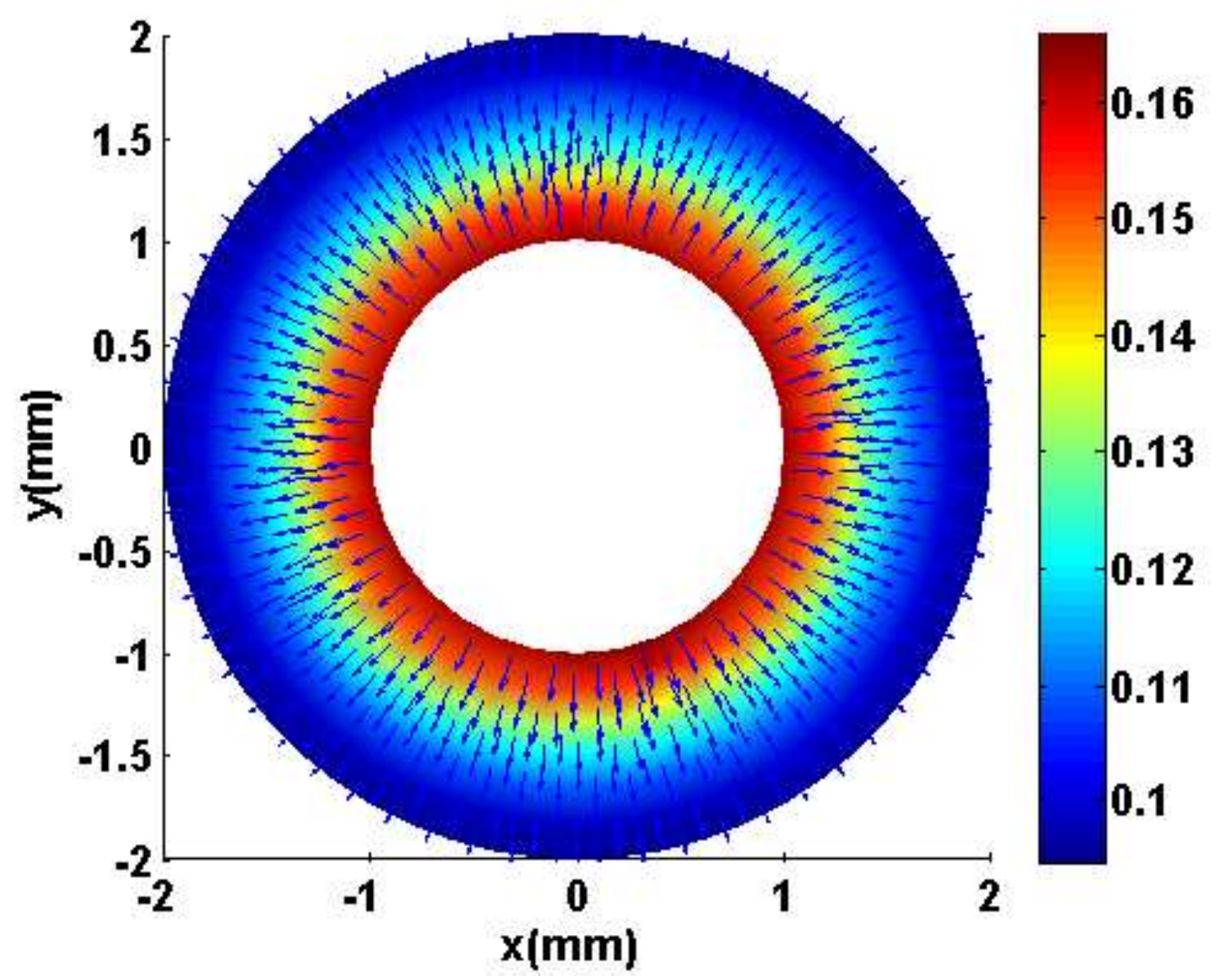}}
   \subfigure[]{
    \label{tecoareal4}
   \includegraphics[width=0.48\columnwidth,draft=false]{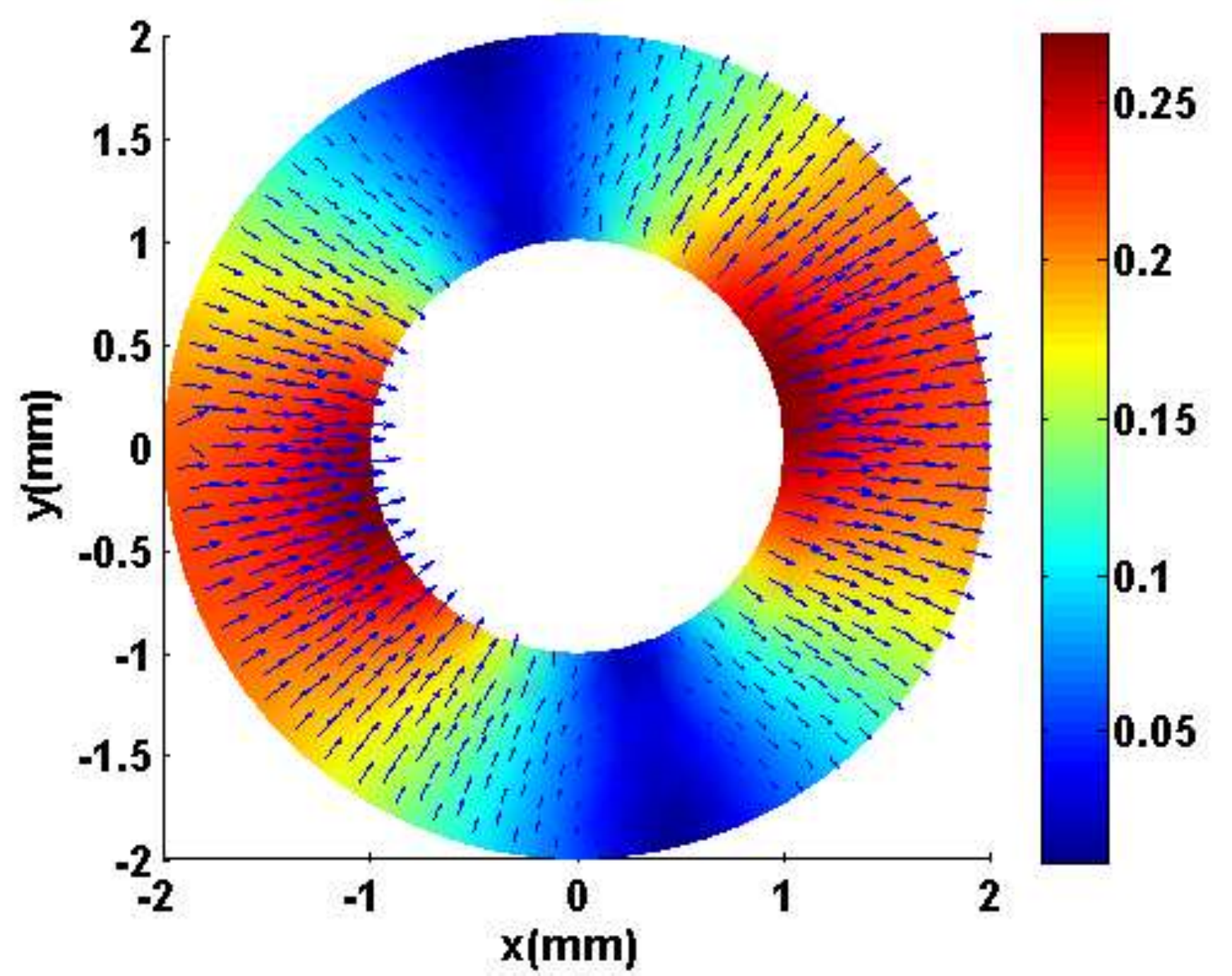}}\\
      \subfigure[]{
    \label{tecoareal5}
   \includegraphics[width=0.48\columnwidth,draft=false]{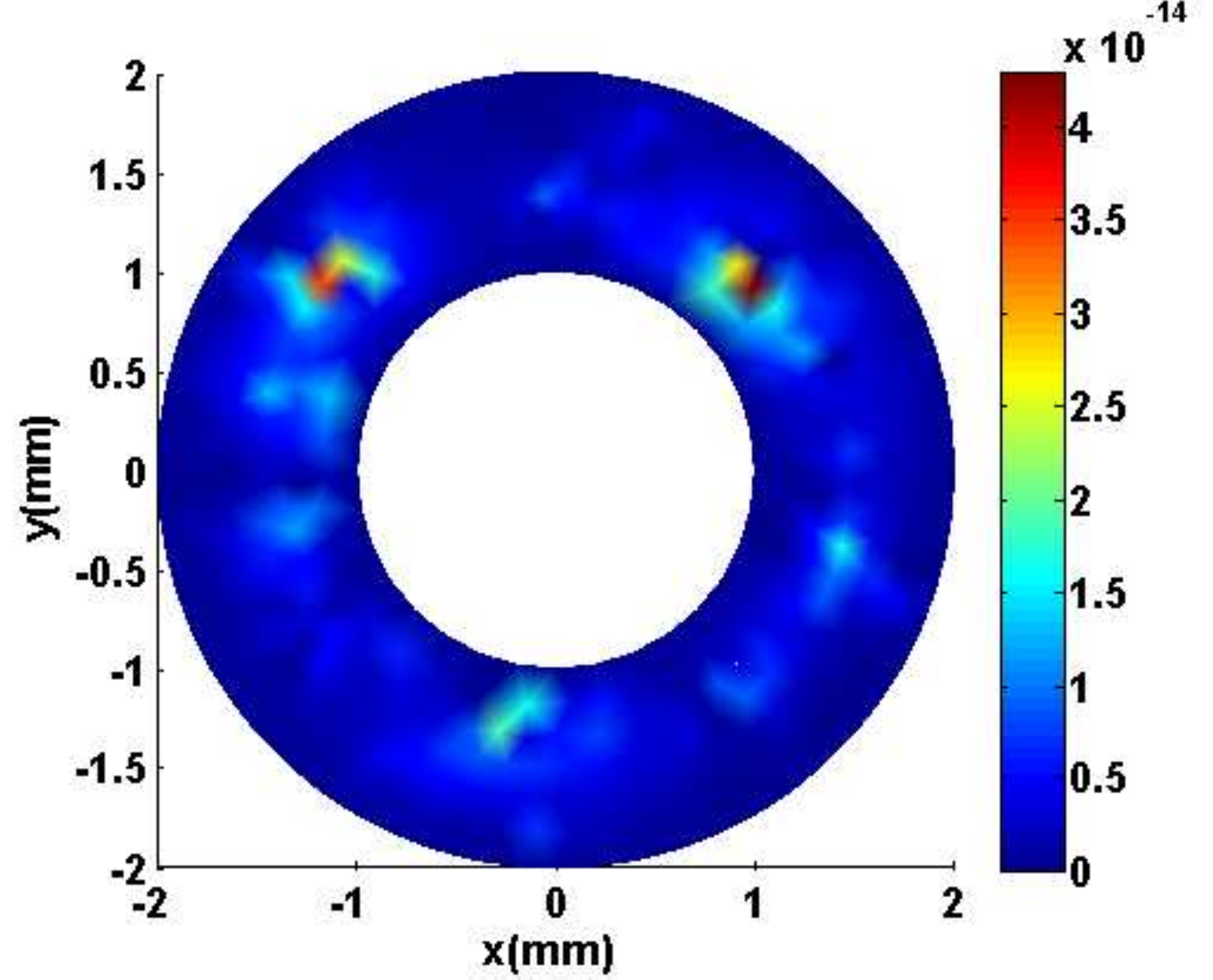}}
   \subfigure[]{
    \label{tecoareal6}
   \includegraphics[width=0.48\columnwidth,draft=false]{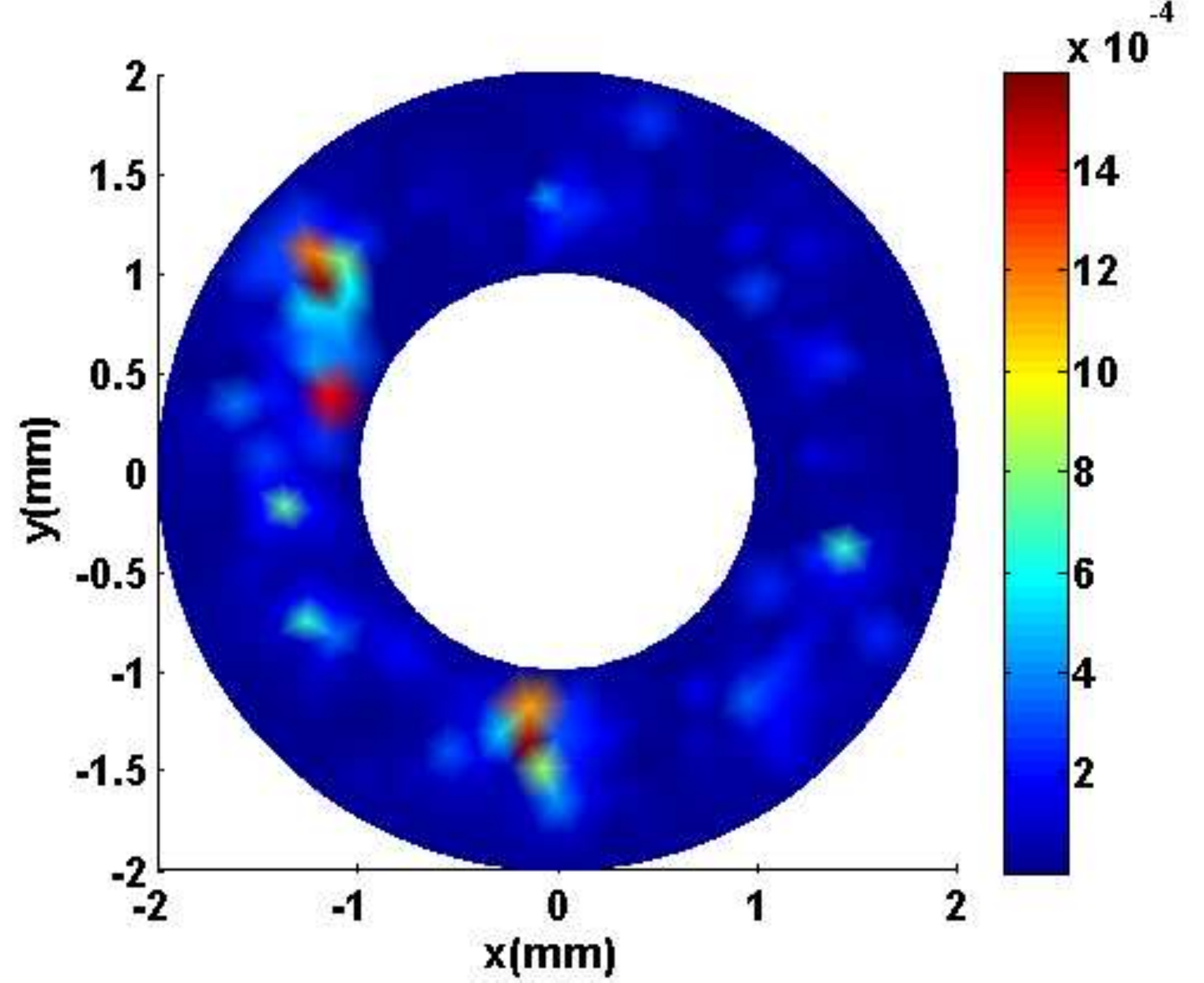}}\\}
 \caption{(a) Magnitude and field distribution of $\mbox{Re}(\ee_{t,1}^{(h)})$.
          (b) Magnitude and field distribution of $\mbox{Re}(\ee_{t,2}^{(h)})$.
          (c) Magnitude and field distribution of $\mbox{Im}(\ee_{t,1}^{(h)})$.
          (d) Magnitude and field distribution of $\mbox{Im}(\ee_{t,2}^{(h)})$.
          (e) Magnitude distribution of $p_{1,h}^{(1)}$.
          (f) Magnitude distribution of $p_{1,h}^{(2)}$.
          These figures are obtained from mixed FEM on the third
          mesh ($h\approx0.2472$).
          The first mode is TEM mode in coaxial waveguide,
          and the second mode is an independent TE mode in coaxial waveguide.}
\label{tecoawaveguide}
\end{figure}

\begin{figure}[ht]
  \centering
  {\subfigure[]{
    \label{tmcoareal1}
   \includegraphics[width=0.48\columnwidth,draft=false]{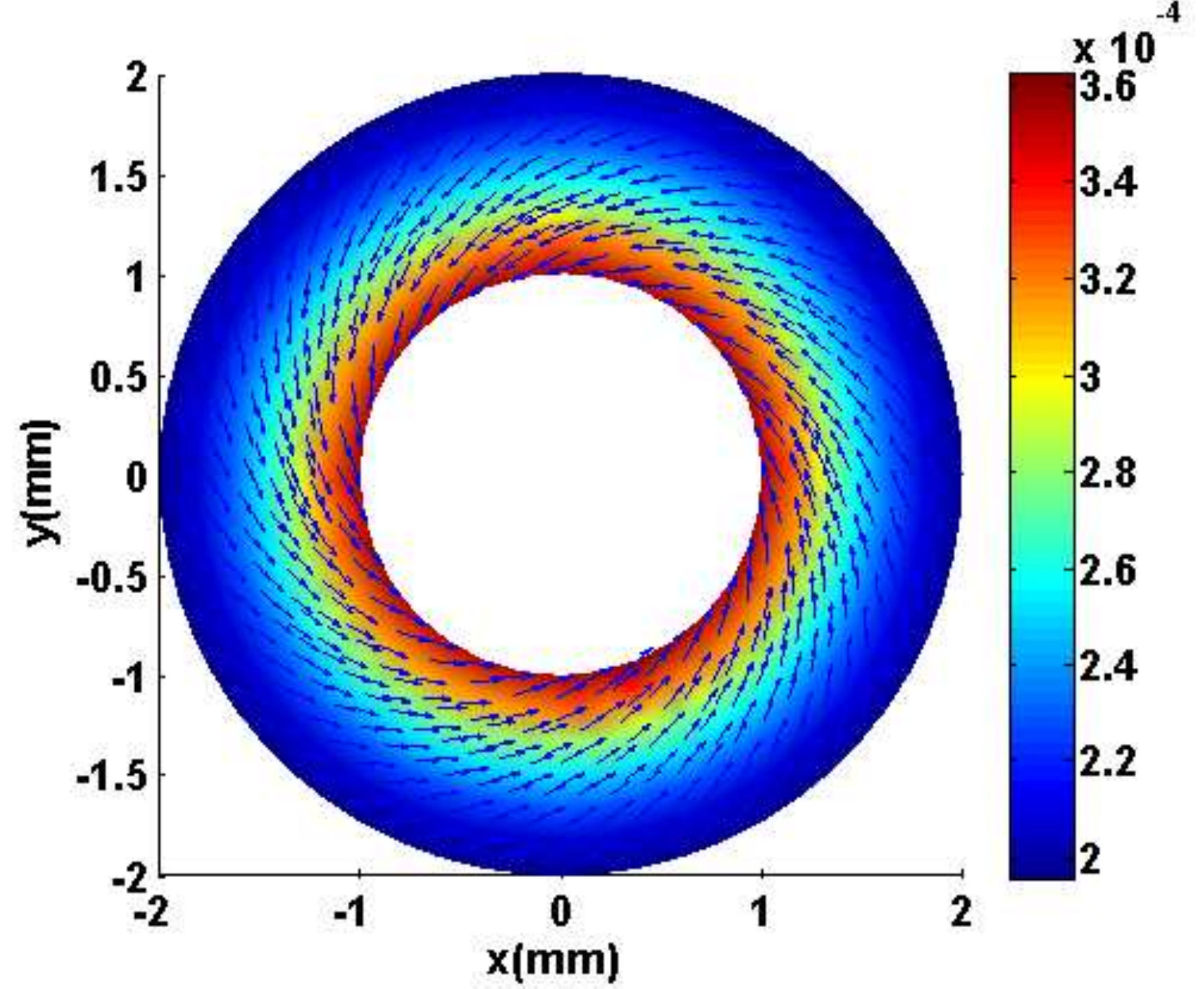}}
     \subfigure[]{
    \label{tmcoareal2}
   \includegraphics[width=0.48\columnwidth,draft=false]{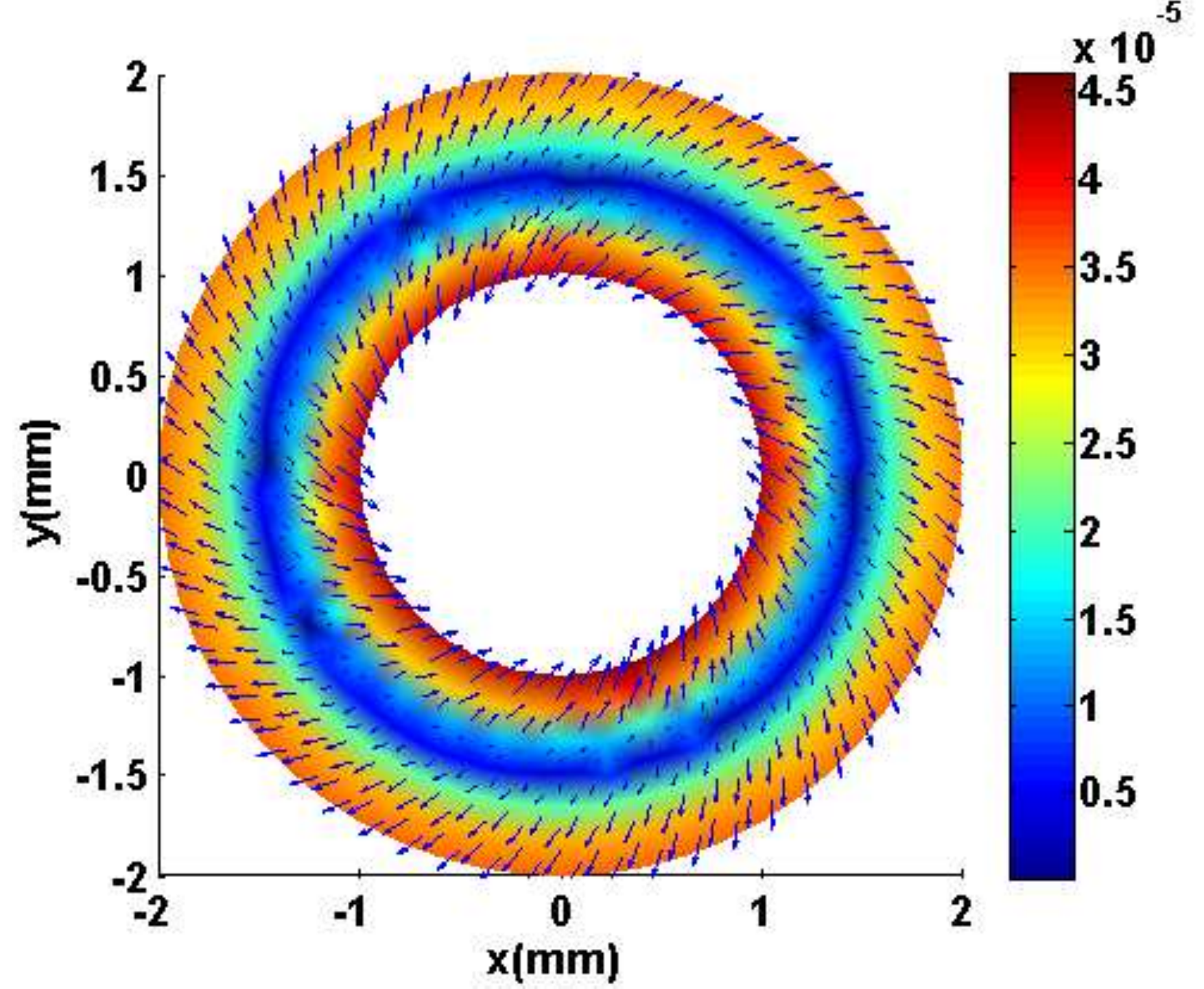}}\\
   \subfigure[]{
    \label{tmcoareal3}
   \includegraphics[width=0.48\columnwidth,draft=false]{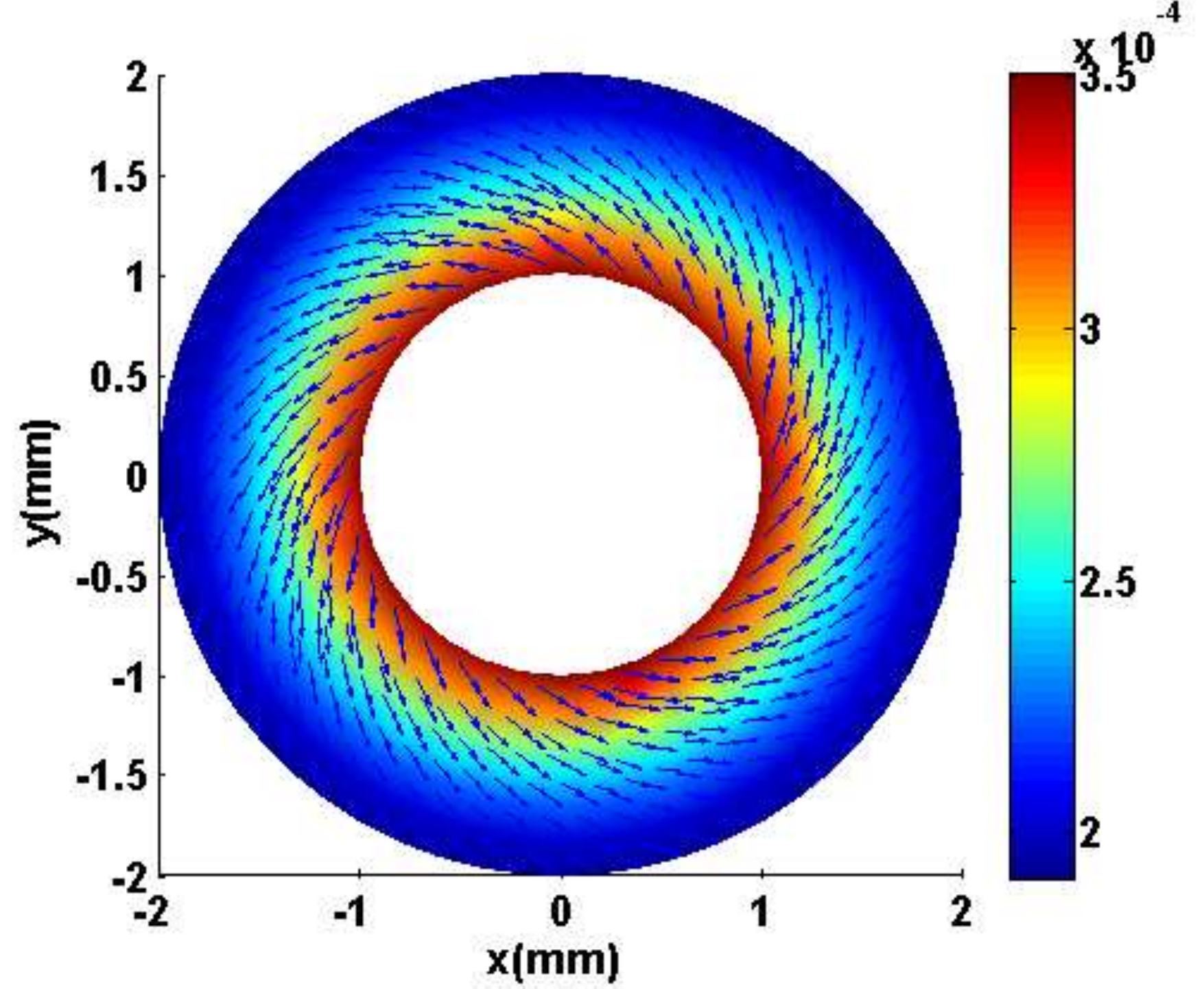}}
   \subfigure[]{
    \label{tmcoareal4}
   \includegraphics[width=0.48\columnwidth,draft=false]{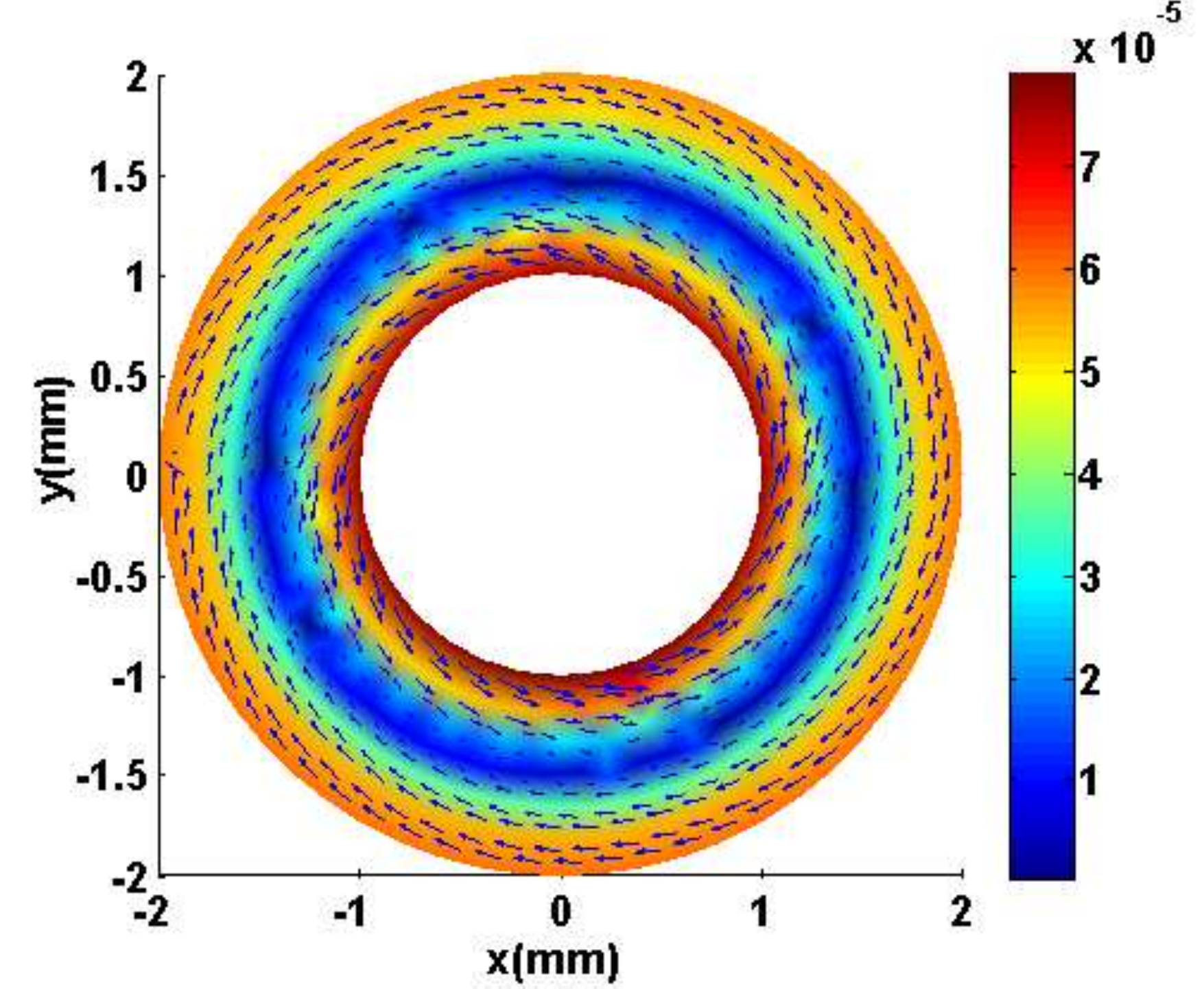}}\\
      \subfigure[]{
    \label{tmcoareal5}
   \includegraphics[width=0.48\columnwidth,draft=false]{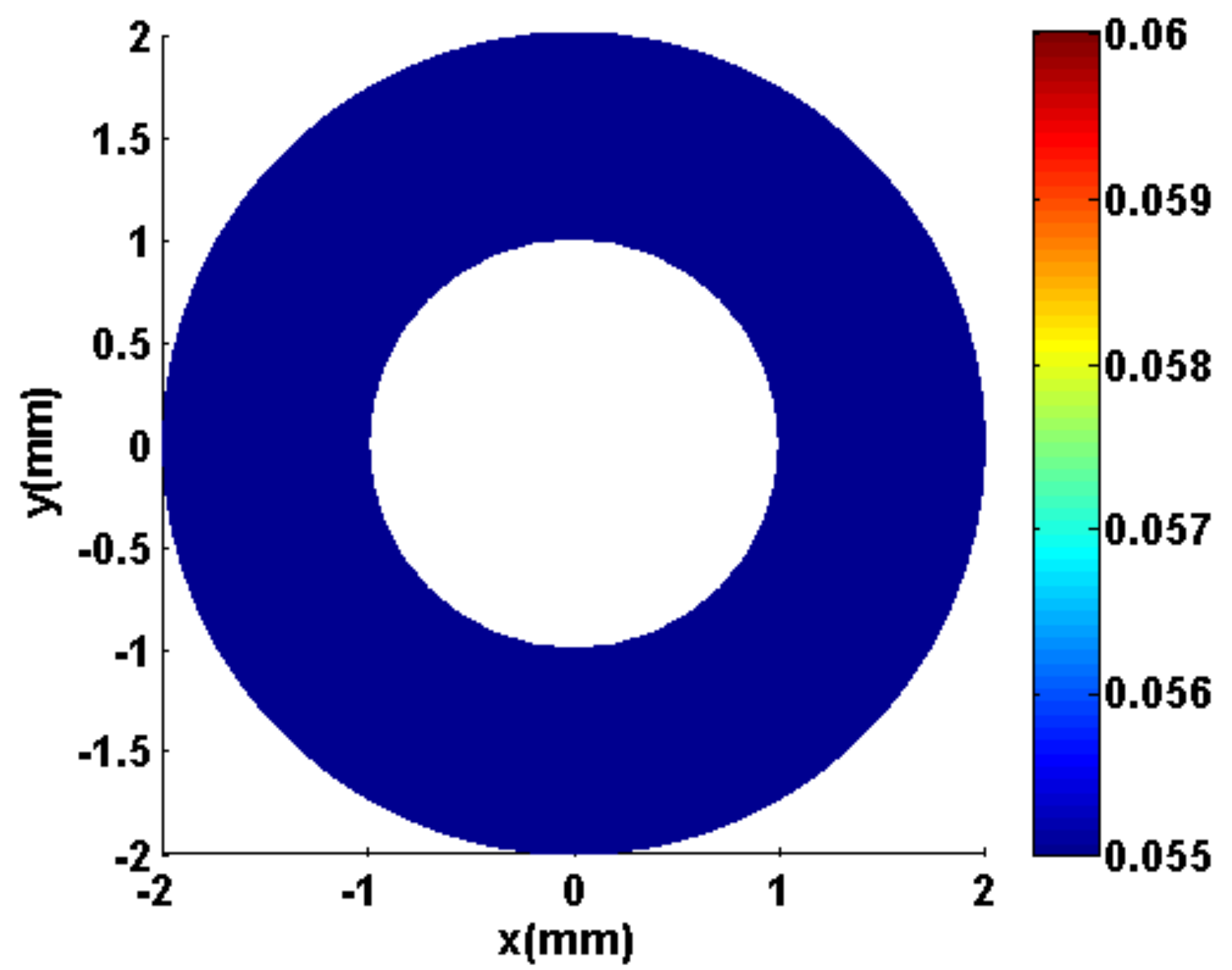}}
   \subfigure[]{
    \label{tmcoareal6}
   \includegraphics[width=0.48\columnwidth,draft=false]{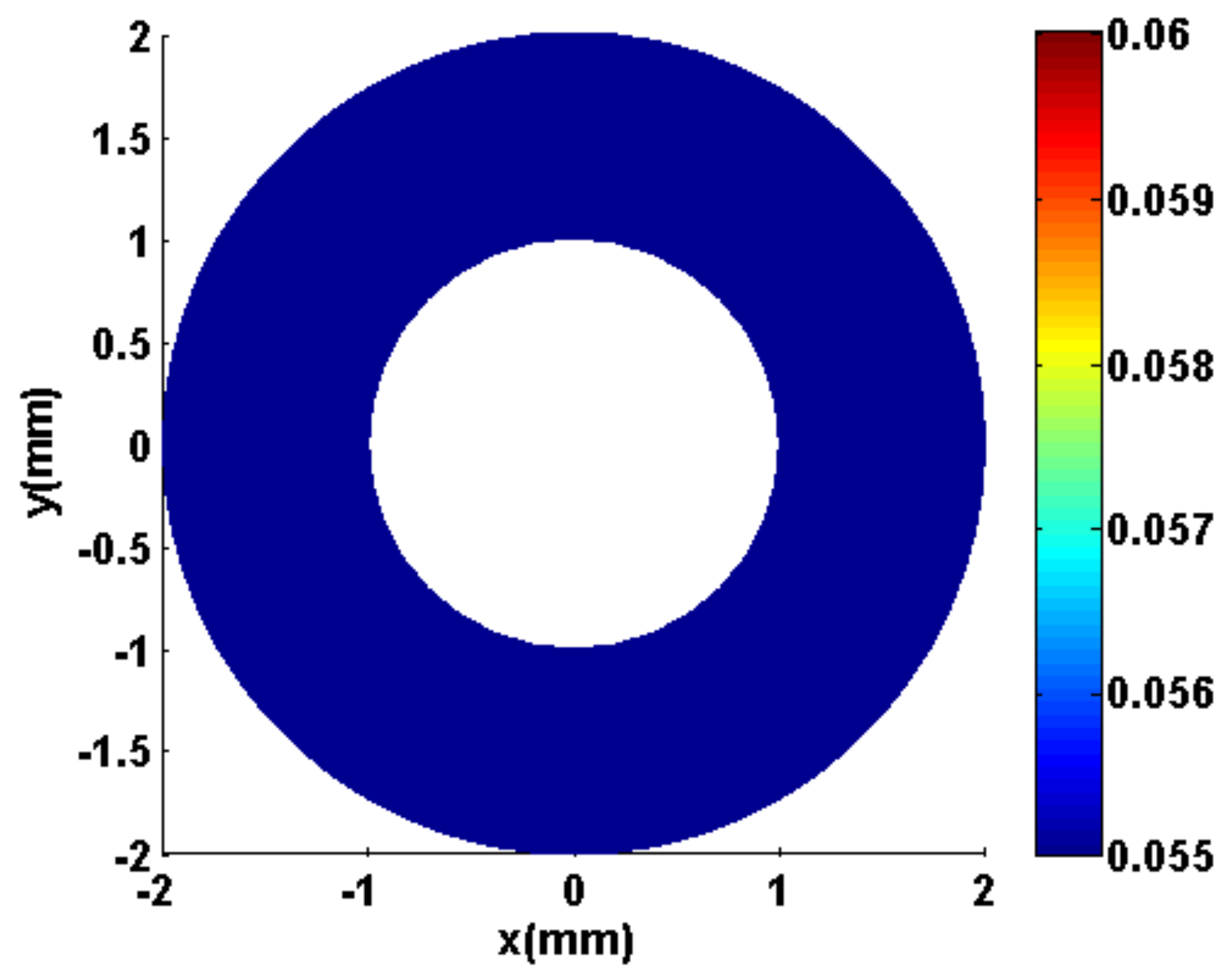}}\\}
 \caption{(a) Magnitude and field distribution of $\mbox{Re}(\h_{t,1}^{(h)})$.
          (b) Magnitude and field distribution of $\mbox{Re}(\h_{t,2}^{(h)})$.
          (c) Magnitude and field distribution of $\mbox{Im}(\h_{t,1}^{(h)})$.
          (d) Magnitude and field distribution of $\mbox{Im}(\h_{t,2}^{(h)})$.
          (e) Magnitude distribution of $p_{2,h}^{(1)}$.
          (f) Magnitude distribution of $p_{2,h}^{(2)}$.
          These figures are obtained from mixed FEM on the third
          mesh ($h\approx0.2472$).
          The first mode is TEM mode in coaxial waveguide,
          and the second mode is an independent TM mode in coaxial waveguide.}
\label{tmcoawaveguide}
\end{figure}

In this section, we shall consider four familiar waveguides, i.e.,  rectangular waveguide, cylindrical waveguide, coaxial waveguide and double-ridge waveguide. These waveguides are filled with homogenous anisotropic lossless medium, which satisfies the above condition \uppercase\expandafter{\romannumeral1} and \uppercase\expandafter{\romannumeral2}. The physical dimensions of these four waveguides are listed in Fig.\ref{figconv}.
Assume that the medium parameters in these four waveguides are
\begin{eqnarray*}
            \d{\ep}=\left[
  \begin{array}{ccc}
    2 &-j &0 \\
    j & 2&0 \\
    0 & 0&1 \\
  \end{array}
\right]\ep_{0},
\quad\d{\mu}=\left[
  \begin{array}{ccc}
     1 &0.5j& 0 \\
     -0.5j&1& 0 \\
    0 & 0&2 \\
  \end{array}
\right]\mu_{0},
\end{eqnarray*}
where $\ep_{0}$ and $\mu_{0}$ are the permittivity and permeability in vacuum, respectively. Obviously, these parameters have satisfied the condition \uppercase\expandafter{\romannumeral1} and \uppercase\expandafter{\romannumeral2}.\\
\indent Let $K_{t,h}^{(q)}, q=0,1,\cdots$ be numerical cut-off wavenumbers obtained from the vector PDEs (\ref{3eqs2}) and (\ref{3eqsm2}). Let $k_{t,h}^{(q)}, q=1,2\cdots$ be numerical wavenumbers obtained from the scalar PDE (\ref{3eqs1}) and (\ref{3eqsm1}).
In the FEM of postprocessing, we need to draw the field distribution in $\Gamma$. For the vector field  $\F=\^x(\F_{x}^{\mbox{real}}+j\F_{x}^{\mbox{imag}})+\^y(\F_{y}^{\mbox{real}}+j\F_{y}^{\mbox{imag}})$,
we denote $\mbox{Re}(\F)=\^x\F_{x}^{\mbox{real}}+\^y\F_{y}^{\mbox{real}}$ and $\mbox{Im}(\F)=\^x\F_{x}^{\mbox{imag}}+\^y\F_{y}^{\mbox{imag}}$. \\
\indent For independent TE modes, we do not list $k_{t,h}^{(0)}$ in Table \uppercase\expandafter{\romannumeral1}-\uppercase\expandafter{\romannumeral4},
because they are not propagable physical modes. For the coaxial waveguide, we find the existence of a TEM mode in Table \uppercase\expandafter{\romannumeral3} and \uppercase\expandafter{\romannumeral7},
because there are two disconnected boundaries in the coaxial waveguide.
Moreover, we have already drawn the magnitude and direction of the electromagnetic field associated with this TEM mode in Fig. \ref{tecoawaveguide} and Fig. \ref{tmcoawaveguide}.
However for other three waveguides, because these three waveguides have connected boundary,
we do not find the existence of any TEM mode, which is in accordance with the waveguide theory.\\
\indent From the Table \uppercase\expandafter{\romannumeral1}-\uppercase\expandafter{\romannumeral4},
we can see that for independent TE modes, all the nonzero numerical cut-off wavenumbers from between the scalar PDE (\ref{3eqs1}) and the vector PDEs (\ref{3eqsm1}) are same roughly, which are agreements with the {\bf {Theorem 1}} in section \uppercase\expandafter{\romannumeral3}.
From the Table \uppercase\expandafter{\romannumeral5}-\uppercase\expandafter{\romannumeral8}, we can see that for independent TM modes, all the nonzero numerical cut-off wavenumbers from between the scalar PDE (\ref{3eqs2}) and the vector PDEs (\ref{3eqsm2}) are same roughly,
which are also agreements with the {\bf {Theorem 2}} in section \uppercase\expandafter{\romannumeral4}. From the Table \uppercase\expandafter{\romannumeral1}-\uppercase\expandafter{\romannumeral8}, we can observe that all the numerical eigenvalues obtained from the scalar PDE (\ref{3eqs1}) and (\ref{3eqsm1}) approximate the exact eigenvalues from above,
which are coincided with the conclusion of min-max principle \cite{strang1973}. However for all numerical eigenvalues obtained from the vector  PDEs (\ref{3eqs2}) and (\ref{3eqsm2}), sometimes they will approximate the exact eigenvalues from above, sometimes approximate the exact eigenvalues from below,
and even swing around from the exact eigenvalues. \\
\indent From the Fig. \ref{terecwaveguide} and Fig. \ref{tecoawaveguide}, we find that the Lagrange multiplier $p_{1,h}^{(l)}\approx0~
(l=1,2)$ corresponding to vector PDEs (\ref{3eqs2}), which is agreement with the theory in section \uppercase\expandafter{\romannumeral5}.  From the Fig. \ref{tmrecwaveguide} and Fig. \ref{tmcoawaveguide}, we can see that the Lagrange multiplier $p_{2,h}^{(l)}=C~
(l=1,2)$ corresponding to vector PDEs (\ref{3eqsm2}), which is also agreement with the theory in section \uppercase\expandafter{\romannumeral5}.\\
\indent At last, from Table \uppercase\expandafter{\romannumeral8}, we see that the numerical accuracy of the first eigenvalue is more worse than one of remaining three eigenvalues,
because the smoothness of the first eigenfunction is very bad \cite{Costabel2000}, this will reduce the rate of convergence for the first numerical eigenvalue.

\begin{table*}[!t]
\renewcommand{\arraystretch}{1.5}
\caption{The first four smallest cut-off wavenumbers ($\times10^3$) from rectangular waveguide (TE modes)}
\centering
\begin{tabular}{ccccccc}
  \hline
 $h(mm)$& 0.195256241897666
 &  0.097628120948833 & 0.048814060474417 & 0.024407030237208& 0.012203515118604
&Trend\\
   \hline
   \hline
 $k_{t,h}^{(1)}$&  1.475789320459881 &1.469813263355673 &  1.468124820384027 &1.467681054232777& 1.467568001103275&$\searrow$\\
 $K_{t,h}^{(1)}$&  1.465403953988726 &1.466972919268568 &  1.467388723990029 &1.467494554209295& 1.467521153945791&$\nearrow$\\
$k_{t,h}^{(2)}$&2.397873950719001  & 2.373928543161163 &  2.367449438270285 & 2.365781632235588& 2.365360394258005&$\searrow$\\
$K_{t,h}^{(2)}$&2.341638022828847  & 2.359141186518130 &  2.363685448494339 & 2.364834868032611& 2.365123211095598&$\nearrow$\\
$k_{t,h}^{(3)}$& 2.448628970844728  & 2.402140321555807 & 2.389806229849646 & 2.386635645572716 &  2.385833987683120&$\searrow$\\
$K_{t,h}^{(3)}$& 2.386120253640742  & 2.385544045675012 & 2.385546681581226 & 2.385559751043563 &  2.385563986757260&$\mbox{Swing}$\\
$k_{t,h}^{(4)}$&3.785795418686821 &  3.640689044398126 & 3.602109334142066& 3.592177442067254 &  3.589665639359895&$\searrow$\\
$K_{t,h}^{(4)}$&3.584051491410978 &  3.587668141693166 & 3.588501184721795& 3.588740989770616 &  3.588803416655484&$\nearrow$\\
 \hline
\end{tabular}
\end{table*}

\begin{table*}[!t]
\renewcommand{\arraystretch}{1.5}
\caption{The first four smallest cut-off wavenumbers ($\times10^3$) from cylindrical waveguide (TE modes)}
\centering
\begin{tabular}{ccccccc}
  \hline
 $h(mm)$& 0.435162097754408
 & 0.318120256910917 & 0.252501014511383 & 0.171260830755769& 0.099425990995704
&Trend\\
   \hline
   \hline
 $k_{t,h}^{(1)}$&  0.481782906245642 &0.480595989735284 &  0.480244222516730 &0.479917287787030& 0.479716334872045&$\searrow$\\
 $K_{t,h}^{(1)}$&  0.482062240494314 &0.480731447714930 &  0.480309890299193 &0.479951716681064& 0.479726122948813&$\searrow$\\
$k_{t,h}^{(2)}$&0.770681543496923  & 0.767985291525049 &  0.767152823640153 & 0.766408960857703& 0.765944500371049&$\searrow$\\
$K_{t,h}^{(2)}$&0.765895295236251  & 0.765781093261658 &  0.765779065046927 & 0.765761030659999& 0.765759132652098&$\searrow$\\
$k_{t,h}^{(3)}$& 0.822512866800075  & 0.817468768730999 & 0.815715297831422 & 0.814371630445448 &  0.813497973724559&$\searrow$\\
$K_{t,h}^{(3)}$& 0.818575391775589  & 0.815610431214487 & 0.814645768992308 & 0.813833505067192 &  0.813347046777445&$\searrow$\\
$k_{t,h}^{(4)}$&1.161993527368619 &  1.147473267954206 & 1.142076505638176& 1.138145314028453 &  1.135646092854336&$\searrow$\\
$K_{t,h}^{(4)}$&1.143769145276075 &  1.138652127539504  & 1.137196475847020& 1.135789867445030 &  1.134985931812746&$\searrow$\\
 \hline
\end{tabular}
\end{table*}

\begin{table*}[!t]
\renewcommand{\arraystretch}{1.5}
\caption{The first five smallest cut-off wavenumbers $(\times10^3)$ from coaxial waveguide (TE modes)}
\centering
\begin{tabular}{ccccccc}
  \hline
 $h(mm)$& 0.472473182342669
 & 0.272133841282597 & 0.247159259110537 & 0.182248215374031& 0.094462741227037
&Trend\\
   \hline
   \hline
 $K_{t,h}^{(0)}$&  0.002061021497E-6 &0.007755090243E-6 &  0.056866086854E-6 &0.038488608349E-6& 0.219807821098E-6&\\
 $k_{t,h}^{(1)}$&  0.403129504958078 &0.401325068249122 &  0.400816169469970 &0.400378940177273 & 0.400081247548716&$\searrow$\\
 $K_{t,h}^{(1)}$&  0.401372769497186 &0.400550372965758 &  0.400264640506610 &0.400108011038728& 0.400004868788086&$\searrow$\\
$k_{t,h}^{(2)}$&0.446491208788408  & 0.437800832627968 &  0.436231024847422 & 0.433984206770929& 0.432413954850601&$\searrow$\\
$K_{t,h}^{(2)}$&0.430961475182653  & 0.431625902509177 &  0.431430711927798 & 0.431622918576020& 0.431713632693326&$\mbox{Swing}$\\
$k_{t,h}^{(3)}$& 0.776672843817230  & 0.772066052975698 & 0.771054257278046 & 0.769924258377010 &  0.769185048375785&$\searrow$\\
$K_{t,h}^{(3)}$& 0.773407091622260  & 0.770960972629768 & 0.769872783975640 & 0.769340846993346 &  0.769010111595282&$\searrow$\\
$k_{t,h}^{(4)}$&0.941131422630201 &  0.911604386951099 & 0.906615497277321& 0.898640732775632 &  0.892784182029530&$\searrow$\\
$K_{t,h}^{(4)}$&0.885483351439040 &  0.888992418911532 & 0.888435742503244& 0.889422642562237 &  0.890035675537056&$\nearrow$\\
 \hline
\end{tabular}
\end{table*}

\begin{table*}[!t]
\renewcommand{\arraystretch}{1.5}
\caption{The first four smallest cut-off wavenumbers ($\times10^3$) from double-ridge waveguide (TE modes)}
\centering
\begin{tabular}{ccccccc}
  \hline
 $h(mm)$& 0.499157708797601
 & 0.307034829640539 & 0.236552539188960 & 0.177102027619670& 0.089101194801190
&Trend\\
   \hline
   \hline
 $k_{t,h}^{(1)}$&  0.298473083901503 &0.296837798083560 &  0.295862095872810 &0.294571635290433& 0.293746949135949&$\searrow$\\
 $K_{t,h}^{(1)}$&  0.289330352346550 &0.290215795563825 &  0.290873024389050 &0.291797090232078& 0.292341396295307&$\nearrow$\\
$k_{t,h}^{(2)}$&0.659141343015736  & 0.653344598256282 &  0.650866137293736 & 0.647691545816427& 0.645837969377189&$\searrow$\\
$K_{t,h}^{(2)}$&0.636364866469366  & 0.638327377661742 &  0.639625708785915 & 0.641630370940502& 0.642778844651659&$\nearrow$\\
$k_{t,h}^{(3)}$& 0.697195503804344  & 0.690640497552013 & 0.687486432230491 & 0.683466011354168 &  0.681128180642666&$\searrow$\\
$K_{t,h}^{(3)}$& 0.671477244947293  & 0.673235103178177 & 0.674510021173172 & 0.676528058790026 &  0.677698537012275&$\nearrow$\\
$k_{t,h}^{(4)}$&0.880308418594263 &  0.861309667546236 & 0.853382847280550& 0.843422634539002 &  0.837647090766269&$\searrow$\\
$K_{t,h}^{(4)}$&0.811516686813181 &  0.817675126146741 & 0.821198506557971& 0.826274208071270 &  0.829194252542041&$\nearrow$\\
 \hline
\end{tabular}
\end{table*}

\begin{table*}[!t]
\renewcommand{\arraystretch}{1.5}
\caption{The first four smallest cut-off wavenumbers ($\times10^3$) from rectangular waveguide (TM modes)}
\centering
\begin{tabular}{ccccccc}
  \hline
 $h(mm)$& 0.195256241897666
 &  0.097628120948833 & 0.048814060474417 & 0.024407030237208& 0.012203515118604
&Trend\\
   \hline
   \hline
 $k_{t,h}^{(1)}$&  5.894531604477931 & 5.811189310908255 &  5.790303984646993 &5.785079349732899& 5.783772980059854&$\searrow$\\
 $K_{t,h}^{(1)}$&  5.751753925681920 &5.773972418995951 &  5.780865626159621 &5.782709191382921& 5.783179649688024&$\nearrow$\\
$k_{t,h}^{(2)}$&9.004880315505943  & 8.729206088628956 &  8.658945775718696 & 8.641300878047364& 8.636884719323140&$\searrow$\\
$K_{t,h}^{(2)}$&8.504444706274022  & 8.599897601466521 &  8.626242822987356 & 8.633095522424496& 8.634831174098622&$\nearrow$\\
$k_{t,h}^{(3)}$& 10.030601701890467  & 9.726233748432023 & 9.651178105103899 & 9.632475354672261 &  9.627803368761551&$\searrow$\\
$K_{t,h}^{(3)}$& 9.390284078717151  & 9.559673599390845 & 9.608995313313448 & 9.621889373857753 &  9.625154021388045&$\nearrow$\\
$k_{t,h}^{(4)}$&12.354007999108354 &  11.782007070090952 & 11.621934827419242& 11.580580289463066 &  11.570156970794526&$\searrow$\\
$K_{t,h}^{(4)}$&11.230485752410944 &  11.486884922583393 & 11.546615650894150& 11.561631137954999 &  11.565411091773859&$\nearrow$\\
 \hline
\end{tabular}
\end{table*}

\begin{table*}[!t]
\renewcommand{\arraystretch}{1.5}
\caption{The first four smallest cut-off wavenumbers ($\times10^3$) from cylindrical waveguide (TM modes)}
\centering
\begin{tabular}{ccccccc}
  \hline
 $h(mm)$& 0.435162097754408
 &  0.318120256910917 & 0.252501014511383 & 0.171260830755769 & 0.099425990995704
&Trend\\
   \hline
   \hline
 $k_{t,h}^{(1)}$&  1.710892465536958 & 1.705229937117969 &  1.703365695670425 &1.701833508509492& 1.700872222278634&$\searrow$\\
 $K_{t,h}^{(1)}$&  1.708130151669107 &1.703909562904990 &  1.702624812725415 &1.701504085656535& 1.700774856957166&$\searrow$\\
$k_{t,h}^{(2)}$&2.748811562159025  & 2.728435217867374 &  2.720656958952525 & 2.714843534766609 & 2.710978247839135&$\searrow$\\
$K_{t,h}^{(2)}$&2.699853520146608  & 2.703738916945631 &  2.706951975137004 & 2.708005426994477& 2.709069682599943&$\nearrow$\\
$k_{t,h}^{(3)}$& 2.750622073669899  & 2.728523176906105 & 2.720979464827615 & 2.714918074075618 &  2.711000912270034&$\searrow$\\
$K_{t,h}^{(3)}$& 2.718342100556592  & 2.713029005422347 & 2.712109453354489 & 2.710766843202556 &  2.709828272360276&$\searrow$\\
$k_{t,h}^{(4)}$&3.723569045100281 &  3.674294064980397 & 3.657351303151986& 3.643844584845230 &  3.635049201403289&$\searrow$\\
$K_{t,h}^{(4)}$&3.592595097303798 &  3.614496263471669 & 3.620490927803921& 3.625553422179841 &  3.629571468542761&$\nearrow$\\
 \hline
\end{tabular}
\end{table*}

\begin{table*}[!t]
\renewcommand{\arraystretch}{1.5}
\caption{The first five smallest cut-off wavenumbers ($\times10^3$) from coaxial waveguide (TM modes)}
\centering
\begin{tabular}{ccccccc}
  \hline
 $h(mm)$& 0.472473182342669
 &  0.272133841282597 & 0.247159259110537 & 0.182248215374031&  0.094462741227037
&Trend\\
   \hline
   \hline
    $K_{t,h}^{(0)}$&  0.066823027505E-5 &0.152203527001E-5 &  0.015769906706E-5 &0.135502281830E-5& 0.24917040751E-5&\\
 $k_{t,h}^{(1)}$&  4.603811174450504 & 4.478909196644620 &  4.462658754701031 &4.436325366189580& 4.423018842494161&$\searrow$\\
 $K_{t,h}^{(1)}$&  4.283083207276171 &4.366132961225913 &  4.382042512252821 &4.400435925712863& 4.411261570362339&$\nearrow$\\
$k_{t,h}^{(2)}$&4.717696213822837  & 4.589919853968560 &  4.570897568984140 & 4.542803941715887& 4.527725382867191&$\searrow$\\
$K_{t,h}^{(2)}$&4.359880505015513  & 4.456850072068127 &  4.472400883228711 & 4.496543178474499& 4.513097773992790&$\nearrow$\\
$k_{t,h}^{(3)}$& 4.726897230985268  & 4.590115719774102 & 4.572194855054836 & 4.542947671779376 &  4.527750693965870&$\searrow$\\
$K_{t,h}^{(3)}$& 4.407953587492303  & 4.477693637268233 & 4.495981667100307 & 4.509426213473715 &  4.516647682087900&$\nearrow$\\
$k_{t,h}^{(4)}$&5.059753397073766 &  4.907645237816515 & 4.882409421458146& 4.847643696333221 &  4.827258453668454&$\searrow$\\
$K_{t,h}^{(4)}$&4.636866073232428 &  4.740140439013424 & 4.758367873685201& 4.786332963346268 &  4.808374979489291&$\nearrow$\\
 \hline
\end{tabular}
\end{table*}

\begin{table*}[!t]
\renewcommand{\arraystretch}{1.5}
\caption{The first four smallest cut-off wavenumbers ($\times10^3$) from double-ridge waveguide (TM modes)}
\centering
\begin{tabular}{ccccccc}
  \hline
 $h(mm)$& 0.195256241897666
 &  0.097628120948833 & 0.048814060474417 & 0.024407030237208& 0.012203515118604
&Trend\\
   \hline
   \hline
 $k_{t,h}^{(1)}$&  4.203895195855380 & 4.124680743300559 &  4.087243448317585 &4.037359496833891& 4.012831978013024&$\searrow$\\
 $K_{t,h}^{(1)}$&  3.783022619082375 &3.859652883157402 &   3.893988098992527 &3.938979866276001& 3.966272994063384&$\nearrow$\\
$k_{t,h}^{(2)}$&4.261953242133240  & 4.180423557759431 &  4.140888274295089 & 4.088927252762932& 4.063393395276099&$\searrow$\\
$K_{t,h}^{(2)}$&3.828157576455789  & 3.904862838812663 &  3.940180926106685 & 3.986480932940768& 4.014931520439401&$\nearrow$\\
$k_{t,h}^{(3)}$& 5.107903446326329  & 4.986683967305610 & 4.955307368211642 & 4.915935482079309 &  4.898370772606942&$\searrow$\\
$K_{t,h}^{(3)}$& 4.759745722763722  & 4.833111151668725  & 4.852557194811825 & 4.874310742059493 &  4.882941715423586&$\nearrow$\\
$k_{t,h}^{(4)}$&5.537715161259985 &  5.435860912946361 & 5.384191950796905& 5.323044003896990 &  5.297012373265178&$\searrow$\\
$K_{t,h}^{(4)}$&5.083946856603320 &  5.163107356166168 & 5.199551654671399& 5.2401919898334819 &  5.261771962435678&$\nearrow$\\
 \hline
\end{tabular}
\end{table*}
\section{Conclusion}
This paper gives a sufficient condition of having independent TE and TM modes in a waveguide filled with homogenous anisotropic lossless medium. In the future, we would like to give a necessary condition of having independent TE and TM modes in a waveguide filled with homogenous anisotropic lossless medium. In that time, the problem that whether there exist independent TE and TM modes in a waveguide filled with homogenous anisotropic lossless medium will be solved thoroughly.

\appendices
\section{Proof of Lemma 1 and Lemma 2}
Proof of {\bf Lemma 1}:
\begin{eqnarray*}
&~~&\curlt\bigg(\mu_{zz}^{-1}\curlt\ee_t\bigg)-k_{t}^2\d{\mu}_{t}^{-1}\ee_t\\
  &=&\frac{j\omega\mu_{zz}^{-1}}{k_{t}^2}\curlt\bigg(\curlt\big(\^z\times(\d{\mu}_{t}\nabla_{t}h_{z})
  \big)\bigg)\\
  &~~&-j\omega\d{\mu}_{t}^{-1}\^z\times(\d{\mu}_{t}\nabla_{t}h_{z})\\
  &=&\frac{j\omega\mu_{zz}^{-1}}{k_{t}^2}\curlt\bigg(\curlt\big(\^z\times(\d{\mu}_{t}\nabla_{t}h_{z})
  \big)\bigg)\\
  &~~&-j\omega\d{\mu}_{t}^{-1}\d{\mu}_{t}\^z\times\nabla_{t}h_{z}\\
  &=&\frac{j\omega\mu_{zz}^{-1}}{k_{t}^2}\curlt\bigg(\curlt\big(\^z\times(\d{\mu}_{t}\nabla_{t}h_{z})
  \big)\bigg)-j\omega\^z\times\nabla_{t}h_{z}\\
  &=&\frac{j\omega\mu_{zz}^{-1}}{k_{t}^2}\curlt\bigg(\^z\cdot(\divt\d{\mu}_{t}\nabla_{t}h_{z})
  \bigg)-j\omega\^z\times\nabla_{t}h_{z}\\
  &=&\frac{j\omega\mu_{zz}^{-1}}{k_{t}^2}\^z\times\bigg(-\gradt(\divt\d{\mu}_{t}\nabla_{t}h_{z})
  \bigg)-j\omega\^z\times\nabla_{t}h_{z}\\
  &=&-\frac{j\omega\mu_{zz}^{-1}}{k_{t}^2}\^z\times\gradt\bigg(\divt(\d{\mu}_{t}\nabla_{t}h_{z})
  +k_{t}^2\mu_{zz}h_{z}\bigg)\\
  &=&0.
\end{eqnarray*}
This has validated the first equation of PDEs (\ref{3eqs2}).\\
\indent Secondly let us verify the correctness of the second equation in PDEs (\ref{3eqs2}).
\begin{eqnarray*}
  \divt\big(\d{\mu}_{t}^{-1}\ee_t\big)&=&\divt\big(\d{\mu}_{t}^{-1}\frac{j\omega}{k_{t}^2}\^z\times(\d{\mu}_{t}\nabla_{t}h_{z})\big)\\
 &=&\frac{j\omega}{k_{t}^2}\divt\big(\d{\mu}_{t}^{-1}\^z\times(\d{\mu}_{t}\nabla_{t}h_{z})\big)\\
 &=&\frac{j\omega}{k_{t}^2}\divt\big(\d{\mu}_{t}^{-1}\d{\mu}_{t}(\^z\times\nabla_{t}h_{z})\big)\\
 &=&\frac{j\omega}{k_{t}^2}\divt\big(\^z\times\nabla_{t}h_{z}\big)\\
 &=&-\frac{j\omega}{k_{t}^2}\bigg(\^z\cdot(\curlt\gradt h_z)\bigg)\\
 &=&0,
\end{eqnarray*}
which validates the second equation of PDEs (\ref{3eqs2}). \\
\indent Finally we validate the boundary condition in PDEs (\ref{3eqs2}).
\begin{eqnarray*}
  \^n\times\ee_t&=&\^n\times\big(\frac{j\omega}{k_{t}^2}\^z\times(\d{\mu}_{t}\nabla_{t}h_{z})\big)\\
   &=&\frac{j\omega}{k_{t}^2}\^n\times\big(\^z\times(\d{\mu}_{t}\nabla_{t}h_{z})\big) \\
   &=&\frac{j\omega}{k_{t}^2}\bigg(\^z\big(\^n\cdot(\d{\mu}_{t}\nabla_{t}h_{z})\big)-\d{\mu}_{t}\nabla_{t}h_{z}(\^n\cdot\^z)\bigg)\\
   &=&0,
\end{eqnarray*}
which validates boundary condition of PDEs (\ref{3eqs2}). The proof of {\bf Lemma 1} is completed.\\

Proof of {\bf Lemma 2}: Firstly we verify the first equation of PDE (\ref{3eqs1}).  When $\omega\neq0$, we have
\begin{eqnarray*}
&~~&k_{t}^2\mu_{zz}h_z+\divt(\d{\mu}_{t}\nabla_{t}h_{z})\\
  &=&k_{t}^2\mu_{zz}\big(\^z\cdot\frac{j\curlt\ee_t}{\omega\mu_{zz}}\big)+\divt\bigg(\d{\mu}_{t}\nabla_{t}(\^z\cdot\frac{j\curlt\ee_t}{\omega\mu_{zz}})\bigg)\\
  &=&\frac{j}{\omega}\bigg(k_{t}^2\^z\cdot\curlt\ee_t+\divt\big(\d{\mu}_{t}\nabla_{t}(\^z\cdot\mu_{zz}^{-1}\curlt\ee_t)\big)\bigg)\\
  &=&\frac{j}{\omega}\bigg(k_{t}^2\^z\cdot\curlt\ee_t+\divt\big(\d{\mu}_{t}\^z\times(\curlt(\mu_{zz}^{-1}\curlt\ee_{t}))\big)\bigg)\\
  &=&\frac{j}{\omega}\bigg(k_{t}^2\^z\cdot\curlt\ee_t+\divt\big(\d{\mu}_{t}\^z\times(k_{t}^2\d{\mu}_{t}^{-1}\ee_{t})\big)\bigg)\\
  &=&\frac{jk_{t}^2}{\omega}\bigg(\^z\cdot\curlt\ee_t+\divt\big(\d{\mu}_{t}\^z\times(\d{\mu}_{t}^{-1}\ee_{t})\big)\bigg)\\
  &=&\frac{jk_{t}^2}{\omega}\bigg(\^z\cdot\curlt\ee_t+\divt\big(\d{\mu}_{t}\d{\mu}_{t}^{-1}\^z\times\ee_{t}\big)\bigg)\\
  &=&\frac{jk_{t}^2}{\omega}\bigg(\^z\cdot\curlt\ee_t+\divt\big(\^z\times\ee_{t}\big)\bigg)\\
  &=&\frac{jk_{t}^2}{\omega}\big(\^z\cdot\curlt\ee_t-\^z\cdot\curlt\ee_t\big)\\
  &=&0,
\end{eqnarray*}
which has already examined the correctness of the first equation of PDE (\ref{3eqs1}); Secondly we verify the boundary condition in PDE (\ref{3eqs1}),
\begin{eqnarray*}
&~~&\^n\cdot(\d{\mu}_{t}\gradt h_z)=\^n\cdot\big(\d{\mu}_{t}\gradt (\^z\cdot\frac{j\curlt\ee_t}{\omega\mu_{zz}})\big)\\
&=&\frac{j}{\omega}\^n\cdot\big(\d{\mu}_{t}\gradt (\^z\cdot\mu_{zz}^{-1}\curlt\ee_t)\big)\\
&=&\frac{j}{\omega}\^n\cdot\big(\d{\mu}_{t}\^z\times(\curlt(\mu_{zz}^{-1}\curlt\ee_t))\big)\\
&=&\frac{j}{\omega}\^n\cdot\big(\d{\mu}_{t}\^z\times(k_{t}^2\d{\mu}_{t}^{-1}\ee_t)\big)\\
&=&\frac{jk_{t}^2}{\omega}\^n\cdot\big(\d{\mu}_{t}\d{\mu}_{t}^{-1}\^z\times\ee_t\big)\\
&=&\frac{jk_{t}^2}{\omega}\^n\cdot\big(\^z\times\ee_t\big)=-\frac{jk_{t}^2}{\omega}\^z\cdot\big(\^n\times\ee_{t}\big)\\
&=&0,
\end{eqnarray*}
which has already examined the boundary condition in PDE (\ref{3eqs1}). The proof of {\bf Lemma 2} is completed.

\section{Abelian Group}
Let $G$ be a set, and $\circ$ be an abstract operation on $G$, if this operation $\circ$ satisfies the following four properties:
\begin{enumerate}
  \item [1.] For every $a,b\in{G}$, there is an element $c\in{G}$, such that $c=a\circ b$, (closure property of the operation $\circ$);
  \item [2.] For every $a,b,c\in{G}$, $(a\circ b)\circ c=a\circ(b\circ c)$,  (associativity);
  \item [3.] There is an identity element $e$, such that $a\circ e=e\circ a=a$, $\forall~a\in{G}$, (the existence of identity element);
  \item [4.] For every $a\in{G}$, there is an element $b\in{G}$, such that
  $a\circ b=b\circ a=e$, (the existence of inverse element);
\end{enumerate}
then we call $(G,\circ)$ is a group. In addition, if $a\circ b=b\circ a$, $\forall~a,b \in{G}$,
we call $(G,\circ)$ is an Abelian group. For details, please see \cite{rotman2010advanced}.

\section{Important Formulas in This Paper}
\begin{enumerate}
  \item [1]
  \begin{gather}
  \int_{\Gamma}\divt(\d{D}\gradt u)\overline{v}dxdy+\int_{\Gamma}\d{D}\gradt u\cdot \gradt\overline{v}dxdy\nonumber\\
  =\int_{\partial\Gamma}\^n\cdot(\d{D}\gradt u)\overline{v}ds \label{scalargreen1}
  \end{gather}
  \item [2]
  \begin{gather}
  \int_{\Gamma}\divt(\d{D}\F)\cdot\overline{q}dxdy+\int_{\Gamma}\d{D}\F\cdot\gradt\overline{q}dxdy\nonumber\\
  =\int_{\partial\Gamma}\^n\cdot(\d{D}\F)\overline{q}ds\label{scalargreen2}
  \end{gather}
  \item [3]
  \begin{gather}
    \int_{\Gamma}\curlt(\alpha\curlt\F_{1})\cdot\overline{\F_{2}}-
    \int_{\Gamma}(\alpha\curlt\F_{1})\cdot(\curlt\overline{\F_{2}})\nonumber\\
    =\int_{\partial\Gamma}\^n\times(\alpha\curlt\F_{1})\cdot\overline{\F_{2}}ds\nonumber\\
    =-\int_{\partial\Gamma}(\alpha\curlt\F_{1})\cdot(\^n\times\overline{\F_{2}})ds\label{vectorgreen}
  \end{gather}
\end{enumerate}
\section*{Acknowledgment}
This work was supported by the National Natural Science Foundation of China under Grant 11371357, Grant 11101381, and Grant 41390453, and by the Fujian Province Natural Science
Foundation under Grant 2013J05060.

\ifCLASSOPTIONcaptionsoff
  \newpage
\fi

\end{document}

%% file: latex_defs.tex
\def\^{\hat}
\def\~{\tilde}
\def\3h{{3\over 2}}

\def\v#1{{\bf#1}}

\def\E{{\v E}}
\def\ee{{\v e}}
\def\h{{\v h}}

\def\H{{\v H}}

\def\B{{\v B}}

\def\F{{\v F}}

\def\W{{\v W}}
\def\N{{\v N}}

\def\ep{\epsilon}

\def\p{\partial}
\def\pd#1#2{{{\partial #1}\over{\partial #2}}}

\def\div{\nabla\cdot}
\def\gradt{\nabla_{t}}
\def\divt{\nabla_{t}\cdot}
\def\curl{\nabla\times}
\def\grad{\nabla}
\def\curlt{\nabla_{t}\times}

\def\eqn#1$${\eqno{{\rm #1}}$$}

\def\ep{\epsilon}

\def\~{\tilde}

\def\^{\hat}

\def\d#1{\overline{\overline{#1}}}

\def\XXint#1#2#3{{\setbox0=\hbox{$#1{#2#3}{\int}$}
     \vcenter{\hbox{$#2#3$}}\kern-.5\wd0}}


%% file: waveguide.bbl
\begin{thebibliography}{10}
\providecommand{\url}[1]{#1}
\csname url@samestyle\endcsname
\providecommand{\newblock}{\relax}
\providecommand{\bibinfo}[2]{#2}
\providecommand{\BIBentrySTDinterwordspacing}{\spaceskip=0pt\relax}
\providecommand{\BIBentryALTinterwordstretchfactor}{4}
\providecommand{\BIBentryALTinterwordspacing}{\spaceskip=\fontdimen2\font plus
\BIBentryALTinterwordstretchfactor\fontdimen3\font minus
  \fontdimen4\font\relax}
\providecommand{\BIBforeignlanguage}[2]{{%
\expandafter\ifx\csname l@#1\endcsname\relax
\typeout{** WARNING: IEEEtran.bst: No hyphenation pattern has been}%
\typeout{** loaded for the language `#1'. Using the pattern for}%
\typeout{** the default language instead.}%
\else
\language=\csname l@#1\endcsname
\fi
#2}}
\providecommand{\BIBdecl}{\relax}
\BIBdecl

\bibitem{hirokawa1998}
J.~Hirokawa and M.~Ando, ``Single-layer feed waveguide consisting of posts for
  plane tem wave excitation in parallel plates,'' \emph{Antennas and
  Propagation, IEEE Transactions on}, vol.~46, no.~5, pp. 625--630, 1998.

\bibitem{Balanis}
C.~A. Balanis, \emph{Advanced Engineering Electromagnetics}.\hskip 1em plus
  0.5em minus 0.4em\relax New York: Wiley, 1989.

\bibitem{rotman2010advanced}
J.~J. Rotman, \emph{Advanced Modern Algebra}.\hskip 1em plus 0.5em minus
  0.4em\relax Upper Saddle River: Prentice Hall, 2002.

\bibitem{Brezzi1991}
F.~Brezzi and M.~Fortin, \emph{Mixed and Hybrid Finite Element Methods}.\hskip
  1em plus 0.5em minus 0.4em\relax New York: Springer-Verlag, 1991.

\bibitem{nede}
J.~C. N\'ed\'elec, ``Mixed finite elements in $\mathbb{R}^3$,'' \emph{Numer.
  Math}, vol.~35, pp. 315--341, 1980.

\bibitem{davidson2010}
D.~B. Davidson, \emph{{Computational Electromagnetics for RF and Microwave
  Engineering}}.\hskip 1em plus 0.5em minus 0.4em\relax Edinburgh: Cambridge
  University Press, 2010.

\bibitem{Jiang2015}
W.~Jiang, N.~Liu, Y.~Qing, and Q.~H. Liu, ``Mixed finite element method for
  resonant cavity problem with complex geometric topology and anisotropic
  media,'' \emph{IEEE Trans. Magnetics}, 2015.

\bibitem{Chew1990}
W.~C. Chew, \emph{Waves and Fields in Inhomogeneous Media}.\hskip 1em plus
  0.5em minus 0.4em\relax New York: Van Nostrand Reinhold, 1990.

\bibitem{Boffi1}
D.~Boffi, ``{Finite Element Approximation of Eigenvalue Problems},'' \emph{Acta
  Numerica}, vol.~19, pp. 1--120, 2010.

\bibitem{chatelin1983}
F.~Chatelin, \emph{Spectral Approximation of Linear Operators}.\hskip 1em plus
  0.5em minus 0.4em\relax New York: Academic Press, 1983.

\bibitem{jiang2014}
W.~Jiang, N.~Liu, Y.~Tang, and Q.~H. Liu, ``{Mixed finite element method for 2D
  vector Maxwell's eigenvalue problem in anisotropic media},'' \emph{Progress
  In Electromagnetics Research}, vol. 148, pp. 159--170, 2014.

\bibitem{ciar}
P.~G. Ciarlet, \emph{The Finite Element Method for Elliptic Problems}.\hskip
  1em plus 0.5em minus 0.4em\relax Amsterdam: North-Holland, 1978.

\bibitem{strang1973}
G.~Strang and G.~J. Fix, \emph{{An Analysis of the Finite Element
  Method}}.\hskip 1em plus 0.5em minus 0.4em\relax New Jersey: Prentice-Hall
  Englewood Cliffs, 1973.

\bibitem{boffi1999}
D.~Boffi, P.~Fernandes, L.~Gastaldi, and I.~Perugia, ``Computational models of
  electromagnetic resonators: Analysis of edge element approximation,''
  \emph{SIAM Journal on Numerical Analysis}, vol.~36, no.~4, pp. 1264--1290,
  1999.

\bibitem{Kikuchi}
F.~Kikuchi, ``Mixed and penalty formulations for finite element analysis of an
  eigenvalue problem in electromagnetism,'' \emph{Computer Methods in Applied
  Mechanics and Engineering}, vol.~64, pp. 509--521, 1987.

\bibitem{Petermonk}
P.~Monk, \emph{{Finite Element Methods for Maxwell's Equations}}.\hskip 1em
  plus 0.5em minus 0.4em\relax Oxford U.K.: Oxford Univ. Press, 2003.

\bibitem{Hiptmair}
R.~Hiptmair, ``{Finite Elements in Computational Electromagnetism},''
  \emph{Acta Numer.}, vol.~11, pp. 237--339, 2002.

\bibitem{Costabel2000}
M.~Costabel and M.~Dauge, ``Singularities of electromagnetic fields in
  polyhedral domains,'' \emph{Arch. Ration. Mech. Anal.}, vol. 151, p.
  221¨C276, 2000.

\end{thebibliography}
